\documentclass[leqno,12pt]{article}

\usepackage{amssymb}
\usepackage{euscript}
\usepackage[dvips]{graphicx}
\usepackage{flafter}
\usepackage{pstricks}

\evensidemargin\oddsidemargin

\begin{document}

 \baselineskip=18pt \setcounter{page}{1}

\renewcommand{\theequation}{\thesection.\arabic{equation}}
\newtheorem{theorem}{Theorem}[section]
\newtheorem{lemma}[theorem]{Lemma}
\newtheorem{proposition}[theorem]{Proposition}
\newtheorem{corollary}[theorem]{Corollary}
\newtheorem{remark}[theorem]{Remark}
\newtheorem{fact}[theorem]{Fact}
\newtheorem{problem}[theorem]{Problem}

\newcommand{\eqnsection}{
\renewcommand{\theequation}{\thesection.\arabic{equation}}
    \makeatletter
    \csname  @addtoreset\endcsname{equation}{section}
    \makeatother}
\eqnsection

\def\r{{\mathbb R}}
\def\e{{\mathbb E}}
\def\p{{\mathbb P}}
\def\P{{\bf P}}
\def\E{{\bf E}}
\def\Q{{\bf Q}}
\def\z{{\mathbb Z}}
\def\T{{\mathbb T}}
\def\G{{\mathbb G}}

\def\ee{\mathrm{e}}
\def\d{\, \mathrm{d}}
\def\deg{{\rm b}}




\vglue50pt

\centerline{\Large\bf A subdiffusive behaviour of recurrent random walk}

\bigskip

\centerline{\Large\bf in random environment on a regular tree}

\bigskip
\bigskip

\centerline{by}

\bigskip

\centerline{Yueyun Hu $\;$and$\;$ Zhan Shi}

\medskip

\centerline{\it Universit\'e Paris XIII \& Universit\'e Paris VI}

\centerline{\tt This version: March 11, 2006}

\bigskip
\bigskip

{\leftskip=2truecm
\rightskip=2truecm
\baselineskip=15pt
\small

\noindent{\slshape\bfseries Summary.} We are interested in the random walk in random environment on an infinite tree. Lyons and Pemantle~\cite{lyons-pemantle} give a precise recurrence/transience criterion. Our paper focuses on the almost sure asymptotic behaviours of a recurrent random walk $(X_n)$ in random environment on a regular tree, which is closely related to Mandelbrot~\cite{mandelbrot}'s multiplicative cascade. We prove, under some general assumptions upon the distribution of the environment, the existence of a new exponent $\nu\in (0, {1\over 2}]$ such that $\max_{0\le i \le n} |X_i|$ behaves asymptotically like $n^{\nu}$. The value of $\nu$ is explicitly formulated in terms of the distribution of the environment.

\bigskip


\bigskip


} 

\bigskip
\bigskip

\section{Introduction}
   \label{s:intro}

Random walk in random environment (RWRE) is a fundamental object in the study of random phenomena in random media. RWRE on $\z$ exhibits rich regimes in the transient case (Kesten, Kozlov and Spitzer~\cite{kesten-kozlov-spitzer}), as well as a slow logarithmic movement in the recurrent case (Sinai~\cite{sinai}). On $\z^d$ (for $d\ge 2$), the study of RWRE remains a big challenge to mathematicians (Sznitman~\cite{sznitman}, Zeitouni~\cite{zeitouni}). The present paper focuses on RWRE on a regular rooted tree, which can be viewed as an infinite-dimensional RWRE. Our main result reveals a rich regime \`a la Kesten--Kozlov--Spitzer, but this time even in the recurrent case; it also strongly suggests the existence of a slow logarithmic regime \`a la Sinai.

Let $\T$ be a $\deg$-ary tree ($\deg\ge 2$) rooted at $e$. For any vertex $x\in \T \backslash \{ e\}$, let ${\buildrel \leftarrow \over x}$ denote the first vertex on the shortest path from $x$ to the root $e$, and $|x|$ the number of edges on this path (notation: $|e|:= 0$). Thus, each vertex $x\in \T \backslash \{ e\}$ has one parent ${\buildrel \leftarrow \over x}$ and $\deg$ children, whereas the root $e$ has $\deg$ children but no parent. We also write ${\buildrel \Leftarrow \over x}$ for the parent of ${\buildrel \leftarrow \over x}$ (for $x\in \T$ such that $|x|\ge 2$).


Let $\omega:= (\omega(x,y), \, x,y\in \T)$ be a family of non-negative random variables such that $\sum_{y\in \T} \omega(x,y)=1$ for any $x\in \T$. Given a realization of $\omega$, we define a Markov chain $X:= (X_n, \, n\ge 0)$ on $\T$ by $X_0 =e$, and whose transition probabilities are
$$
P_\omega(X_{n+1}= y \, | \, X_n =x) = \omega(x, y) .
$$

\noindent Let $\P$ denote the distribution of $\omega$, and let $\p (\cdot) := \int P_\omega (\cdot) \P(\! \d \omega)$. The process $X$ is a $\T$-valued RWRE. (By informally taking $\deg=1$, $X$ would become a usual RWRE on the half-line $\z_+$.)

For general properties of tree-valued processes, we refer to Peres~\cite{peres} and Lyons and Peres~\cite{lyons-peres}. See also Duquesne and Le Gall~\cite{duquesne-le-gall} and Le Gall~\cite{le-gall} for continuous random trees. For a list of motivations to study RWRE on a tree, see Pemantle and Peres~\cite{pemantle-peres1}, p.~106.

We define
\begin{equation}
    A(x) := {\omega({\buildrel \leftarrow \over x},
    x) \over \omega({\buildrel \leftarrow \over x},
    {\buildrel \Leftarrow \over x})} , \qquad
    x\in \T, \; |x|\ge 2.
    \label{A}
\end{equation}

\noindent Following Lyons and Pemantle~\cite{lyons-pemantle}, we assume throughout the paper that $(\omega(x,\bullet))_{x\in \T\backslash \{ e\} }$ is a family of i.i.d.\ {\it non-degenerate} random vectors and that $(A(x), \; x\in \T, \; |x|\ge 2)$ are identically distributed. We also assume the existence of $\varepsilon_0>0$ such that $\omega(x,y) \ge \varepsilon_0$ if either $x= {\buildrel \leftarrow \over y}$ or $y= {\buildrel \leftarrow \over x}$, and $\omega(x,y) =0$ otherwise; in words, $(X_n)$ is a nearest-neighbour walk, satisfying an ellipticity condition.

Let $A$ denote a generic random variable having the common distribution of $A(x)$ (for $|x| \ge 2$). Define
\begin{equation}
    p := \inf_{t\in [0,1]} \E (A^t) .
    \label{p}
\end{equation}

\noindent We recall a recurrence/transience criterion from Lyons and Pemantle (\cite{lyons-pemantle}, Theorem 1 and Proposition 2).

\bigskip

\noindent {\bf Theorem A (Lyons and Pemantle~\cite{lyons-pemantle})} {\it With $\p$-probability one, the walk $(X_n)$ is recurrent or transient, according to whether $p\le {1\over \deg}$ or $p>{1\over \deg}$. It is, moreover, positive recurrent if $p<{1\over \deg}$.}

\bigskip

We study the recurrent case $p\le {1\over \deg}$ in this paper. Our first result, which is not deep, concerns the positive recurrent case $p< {1\over \deg}$.

\medskip

\begin{theorem}
 \label{t:posrec}
 If $p<{1\over \deg}$, then
 \begin{equation}
     \lim_{n\to \infty} \, {1\over \log n} \,
     \max_{0\le i\le n} |X_i| =
     {1\over \log[1/(q\deg)]},
     \qquad \hbox{\rm $\p$-a.s.},
     \label{posrec}
 \end{equation}
 where the constant $q$ is defined in $(\ref{q})$,
 and lies in $(0, {1\over \deg})$ when $p<{1\over
 \deg}$.
\end{theorem}

\medskip

Despite the warning of Pemantle~\cite{pemantle} (``there are many papers proving results on trees as a somewhat unmotivated alternative \dots to Euclidean space"), it seems to be of particular  interest to study the more delicate situation $p={1\over \deg}$ that turns out to possess rich regimes. We prove that, similarly to the Kesten--Kozlov--Spitzer theorem for {\it transient} RWRE on the line, $(X_n)$ enjoys, even in the recurrent case, an interesting subdiffusive behaviour.

To state our main result, we define
\begin{eqnarray}
    \kappa
 &:=& \inf\left\{ t>1: \; \E(A^t) = {1\over \deg}
    \right\} \in (1, \infty], \qquad (\inf
    \emptyset=\infty)
    \label{kappa}
    \\
    \psi(t)
 &:=& \log \E \left( A^t \right) , \qquad t\ge 0.
    \label{psi}
\end{eqnarray}

\noindent We use the notation $a_n \approx b_n$ to denote $\lim_{n\to \infty} \, {\log a_n \over \log b_n} =1$.

\medskip

\begin{theorem}
 \label{t:nullrec}
 If $p={1\over \deg}$ and if $\psi'(1)<0$, then
 \begin{equation}
     \max_{0\le i\le n} |X_i| \; \approx\; n^\nu,
     \qquad \hbox{\rm $\p$-a.s.},
     \label{nullrec}
 \end{equation}
 where $\nu=\nu(\kappa)$ is defined by
 \begin{equation}
     \nu := 1- {1\over \min\{ \kappa, 2\} } =
     \left\{
     \begin{array}{ll}
         (\kappa-1)/\kappa,
      & \mbox{if $\;\kappa \in (1,2]$},
         \\
         \\
         1/2
      & \mbox{if $\;\kappa \in (2, \infty].$}
     \end{array} \right.
     \label{theta}
 \end{equation}
\end{theorem}

\medskip

\noindent {\bf Remark.} (i) It is known (Menshikov and Petritis~\cite{menshikov-petritis}) that if $p={1\over \deg}$ and $\psi'(1)<0$, then for $\P$-almost all environment $\omega$, $(X_n)$ is
null recurrent.

(ii) For the value of $\kappa$, see Figure 1. Under the assumptions $p={1\over \deg}$ and $\psi'(1)<0$, the value of $\kappa$ lies in $(2, \infty]$ if and only if $\E (A^2) <  {1\over \deg}$; and $\kappa=\infty$ if moreover $\hbox{ess sup}(A) \le 1$.

(iii) Since the walk is recurrent, $\max_{0\le i\le n} |X_i|$ cannot be replaced by $|X_n|$ in (\ref{posrec}) and (\ref{nullrec}).

(iv) Theorem \ref{t:nullrec}, which could be considered as a (weaker) analogue of the Kesten--Kozlov--Spitzer theorem, shows that tree-valued RWRE has even richer regimes than RWRE on $\z$. In fact, recurrent RWRE on $\z$ is of order of magnitude $(\log n)^2$, and has no $n^a$ (for $0<a<1$) regime.

(v) The case $\psi'(1)\ge 0$ leads to a phenomenon similar to Sinai's slow movement, and is studied in a forthcoming paper.

\bigskip


\begin{figure}[h]
\centerline{%
\hbox{%
\psset{unit=4.5mm} \pspicture(0,0)(39,10) \Cartesian(4.5mm,4.5mm)
\psline{->}(2,2)(14,2) \psline{->}(2,2)(2,10)
\psline[linewidth=1pt,linestyle=dashed]{-}(2,4)(14,4)
\psline[linewidth=1pt,linestyle=dashed]{-}(5.6,2)(5.6,4)
\psbezier[linewidth=1.5pt]{-}(2,8)(5,2)(11,2.5)(13, 5)
\rput[rb](1.8,1.5){$0$} \rput[rb](3.5, 10){$\E(A^t)$}
\rput[rb](1.8,7.8){$1$} \rput[rb](1.8,3.5){$1\over \deg$}
\rput[tl](5.4,1.8){$1$}
 \rput[tr](12.3,1.5){$\kappa$}
 \rput[tr](14, 1.5){$t$}
 \psline[linewidth=1pt,linestyle=dashed]{-}(12,2)(12,4)
\psline{->}(22,2)(36,2) \psline{->}(22,2)(22,10)
\psline[linewidth=1pt,linestyle=dashed]{-}(22,4)(36,4)
\psline[linewidth=1pt,linestyle=dashed]{-}(26,2)(26,4)
\psbezier[linewidth=1.5pt]{-}(22,8)(26,2)(30, 2.9 )(35, 2.6)
\rput[rb](21.8,1.5){$0$} \rput[rb](21.8,7.8){$1$}
\rput[rb](23,10){$\E(A^t)$} \rput[rb](21.8,3.5){$1\over \deg$}
\rput[tl](25.8,1.8){$1$} \rput[lt](30,1.5){$\kappa=\infty$}
\rput[tr](36, 1.5){$t$}
\endpspicture}}
\caption{ \label{fig:l-lt-1} \textsf{The function $t \to \E(A^t)$
in the   case $\psi'(1) < 0 $ and $p={1\over\deg}$.}}
\end{figure}

\bigskip

The rest of the paper is organized as follows. Section \ref{s:posrec} is devoted to the proof of Theorem \ref{t:posrec}. In Section \ref{s:proba}, we collect some elementary inequalities, which will be of frequent use later on. Theorem \ref{t:nullrec} is proved in Section \ref{s:nullrec}, by means of a result (Proposition \ref{p:beta-gamma}) concerning the solution of a recurrence equation which is closely related to Mandelbrot's multiplicative cascade. We prove Proposition \ref{p:beta-gamma} in Section \ref{s:beta-gamma}.

Throughout the paper, $c$ (possibly with a subscript) denotes a finite and positive constant; we write $c(\omega)$ instead of $c$ when the value of $c$ depends on the environment $\omega$.

\section{Proof of Theorem \ref{t:posrec}}
\label{s:posrec}

We first introduce the constant $q$ in the statement of Theorem \ref{t:posrec}, which is defined without the assumption $p< {1\over \deg}$. Let
$$
\varrho(r) := \inf_{t\ge 0} \left\{ r^{-t} \, \E(A^t) \right\} , \qquad r>0.
$$

\noindent Let $\underline{r} >0$ be such that
$$
\log \underline{r} = \E(\log A) .
$$

\noindent We mention that $\varrho(r)=1$ for $r\in (0, \underline{r}]$, and that $\varrho(\cdot)$ is continuous and (strictly) decreasing on $[\underline{r}, \, \Theta)$ (where $\Theta:= \hbox{ess sup}(A) < \infty$), and $\varrho(\Theta) = \P (A= \Theta)$. Moreover, $\varrho(r)=0$ for $r> \Theta$. See Chernoff~\cite{chernoff}.

We define
$$
\overline{r} := \inf\left\{ r>0: \; \varrho(r) \le {1\over \deg} \right\}.
$$

\noindent Clearly, $\underline{r} < \overline{r}$.

We define
\begin{equation}
    q:= \sup_{r\in [\underline{r}, \, \overline{r}]}
    r \varrho(r).
    \label{q}
\end{equation}

\noindent The following elementary lemma tells us that, instead of $p$, we can also use $q$ in the recurrence/transience criterion of Lyons and Pemantle.

\medskip

\begin{lemma}
 \label{l:pq}
 We have $q>{1\over \deg}$ $($resp., $q={1\over
 \deg}$, $q<{1\over \deg})$ if and only if
 $p>{1\over \deg}$ $($resp., $p={1\over \deg}$,
 $p<{1\over \deg})$.
\end{lemma}

\medskip

\noindent {\it Proof of Lemma \ref{l:pq}.} By Lyons and Pemantle (\cite{lyons-pemantle}, p.~129), $p= \sup_{r\in (0, \, 1]} r \varrho (r)$. Since $\varrho(r) =1$ for $r\in (0, \, \underline{r}]$, there exists $\min\{\underline{r}, 1\}\le r^* \le 1$ such that $p= r^* \varrho (r^*)$.

(i) Assume $p<{1\over \deg}$. Then $\varrho (1) \le \sup_{r\in (0, \, 1]} r \varrho (r) = p < {1\over \deg}$, which, by definition of $\overline{r}$, implies $\overline{r} < 1$. Therefore, $q \le p <{1\over \deg}$.

(ii) Assume $p\ge {1\over \deg}$. We have $\varrho (r^*) \ge p \ge {1\over \deg}$, which yields $r^* \le \overline{r}$. If $\underline{r} \le 1$, then $r^*\ge \underline{r}$, and thus $p=r^* \varrho (r^*) \le q$. If $\underline{r} > 1$, then $p=1$, and thus $q\ge \underline{r}\, \varrho (\underline{r}) = \underline{r} > 1=p$.

We have therefore proved that $p\ge {1\over \deg}$ implies $q\ge p$.

If moreover $p>{1\over \deg}$, then $q \ge p>{1\over \deg}$.

(iii) Assume $p={1\over \deg}$. We already know from (ii) that $q \ge p$.

On the other hand, $\varrho (1) \le \sup_{r\in (0, \, 1]} r \varrho (r) = p = {1\over \deg}$, implying $\overline{r} \le 1$. Thus $q \le p$.

As a consequence, $q=p={1\over \deg}$.\hfill$\Box$

\bigskip

Having defined $q$, the next step in the proof of Theorem \ref{t:posrec} is to compute invariant measures $\pi$ for $(X_n)$. We first introduce some notation on the tree. For any $m\ge 0$, let
$$
\T_m := \left\{x \in \T: \; |x| = m \right\} .
$$

\noindent For any $x\in \T$, let $\{ x_i \}_{1\le i\le \deg}$ be the set of children of $x$.

If $\pi$ is an invariant measure, then
$$
\pi (x) = {\omega ({\buildrel \leftarrow \over x}, x) \over \omega (x, {\buildrel \leftarrow \over x})} \, \pi({\buildrel \leftarrow \over x}), \qquad \forall \, x\in \T \backslash \{ e\}.
$$

\noindent By induction, this leads to (recalling $A$ from (\ref{A})): for $x\in \T_m$ ($m\ge 1$),
$$
\pi (x) = {\pi(e)\over \omega (x, {\buildrel \leftarrow \over x})} {\omega (e, x^{(1)}) \over A(x^{(1)})}  \exp\left( \, \sum_{z\in ]\! ] e, x]\! ]} \log A(z) \right) ,
$$

\noindent where $]\! ] e, x]\! ]$ denotes the shortest path $x^{(1)}$, $x^{(2)}$, $\cdots$, $x^{(m)} =: x$ from the root $e$ (but excluded) to the vertex $x$. The identity holds for {\it any} choice of $(A(e_i), \, 1\le i\le \deg)$. We choose $(A(e_i), \, 1\le i\le \deg)$ to be a random vector independent of $(\omega(x,y), \, |x|\ge 1, \, y\in \T)$, and distributed as $(A(x_i), \, 1\le i\le \deg)$, for any $x\in \T_m$ with $m\ge 1$.

By the ellipticity condition on the environment, we can take $\pi(e)$ to be sufficiently small so that for some $c_0\in (0, 1]$,
\begin{equation}
    c_0\, \exp\left( \, \sum_{z\in
    ]\! ] e, x]\! ]} \log A(z) \right) \le \pi (x)
    \le \exp\left( \, \sum_{z\in
    ]\! ] e, x]\! ]} \log A(z) \right) .
    \label{pi}
\end{equation}

By Chebyshev's inequality, for any $r>\underline{r}$,
\begin{equation}
    \max_{x\in \T_n} \P \left\{ \pi (x)
    \ge r^n\right\} \le \varrho(r)^n.
    \label{chernoff}
\end{equation}

\noindent Since $\# \T_n = \deg^n$, this gives $\E (\#\{ x\in \T_n: \; \pi (x)\ge r^n \} ) \le \deg^n \varrho(r)^n$. By Chebyshev's inequality and the Borel--Cantelli lemma, for any $r>\underline{r}$ and $\P$-almost surely for all large $n$,
\begin{equation}
    \#\left\{ x\in \T_n: \; \pi (x) \ge r^n \right\}
    \le n^2 \deg^n \varrho(r)^n.
    \label{Jn-ub1}
\end{equation}

On the other hand, by (\ref{chernoff}),
$$
\P \left\{ \exists x\in \T_n: \pi (x) \ge r^n\right\} \le \deg^n \varrho (r)^n.
$$

\noindent For $r> \overline{r}$, the expression on the right-hand side is summable in $n$. By the Borel--Cantelli lemma, for any $r>\overline{r}$ and $\P$-almost surely for all large $n$,
\begin{equation}
    \max_{x\in \T_n} \pi (x) < r^n.
    \label{Jn-ub}
\end{equation}

\bigskip

\noindent {\it Proof of Theorem \ref{t:posrec}: upper bound.} Fix $\varepsilon>0$ such that $q+ 3\varepsilon < {1\over \deg}$.

We follow the strategy given in Liggett (\cite{liggett}, p.~103) by introducing a positive recurrent birth-and-death chain $(\widetilde{X_j}, \, j\ge 0)$, starting from $0$, with transition probability from $i$ to $i+1$ (for $i\ge 1$) equal to
$$
{1\over \widetilde{\pi} (i)} \, \sum_{x\in \T_i} \pi(x) (1- \omega(x, {\buildrel \leftarrow \over x})) ,
$$

\noindent where $\widetilde{\pi} (i) := \sum_{x\in \T_i} \pi(x)$. We note that $\widetilde{\pi}$ is a finite invariant measure for $(\widetilde{X_j})$.
Let
$$
\tau_n := \inf \left\{ i\ge 1: \, X_i \in \T_n\right\}, \qquad n\ge 0.
$$

\noindent By Liggett (\cite{liggett}, Theorem II.6.10), for any $n\ge 1$,
$$
P_\omega (\tau_n< \tau_0) \le \widetilde{P}_\omega (\widetilde{\tau}_n< \widetilde{\tau}_0),
$$

\noindent where $\widetilde{P}_\omega (\widetilde{\tau}_n< \widetilde{\tau}_0)$ is the probability that $(\widetilde{X_j})$ hits $n$ before returning to $0$. According to Hoel et al.\ (\cite{hoel-port-stone}, p.~32, Formula (61)),
$$
\widetilde{P}_\omega (\widetilde{\tau}_n< \widetilde{\tau}_0) = c_1(\omega) \left( \, \sum_{i=0}^{n-1} {1\over \sum_{x\in \T_i} \pi(x) (1- \omega(x, {\buildrel \leftarrow \over x}))}Ê\right)^{\! \! -1} ,
$$

\noindent where $c_1(\omega)\in (0, \infty)$ depends on $\omega$. We arrive at the following estimate: for any $n\ge 1$,
\begin{equation}
    P_\omega (\tau_n< \tau_0) \le c_1(\omega) \,
    \left( \, \sum_{i=0}^{n-1} {1\over \sum_{x\in
    \T_i} \pi(x)}Ê\right)^{\! \! -1} .
    \label{liggett}
\end{equation}

We now estimate $\sum_{i=0}^{n-1} {1\over \sum_{x\in \T_i} \pi(x)}$. For any fixed $0=r_0< \underline{r} < r_1 < \cdots < r_\ell = \overline{r} <r_{\ell +1}$,
$$
\sum_{x\in \T_i} \pi(x) \le \sum_{j=1}^{\ell+1} (r_j)^i \# \left\{ x\in \T_i:  \pi(x) \ge (r_{j-1})^i \right\} + \sum_{x\in \T_i: \, \pi(x) \ge (r_{\ell +1})^i} \pi(x).
$$

\noindent By (\ref{Jn-ub}), $\sum_{x\in \T_i: \, \pi(x) \ge (r_{\ell +1})^i} \pi(x) =0$ $\P$-almost surely for all large $i$. It follows from (\ref{Jn-ub1}) that $\P$-almost surely, for all large $i$,
$$
\sum_{x\in \T_i} \pi(x) \le (r_1)^i \deg^i +  \sum_{j=2}^{\ell+1} (r_j)^i i^2 \, \deg^i \varrho (r_{j-1})^i.
$$

\noindent Recall that $q= \sup_{r\in [\underline{r}, \, \overline{r}] } r \, \varrho(r) \ge \underline{r} \, \varrho (\underline{r}) = \underline{r}$. We choose $r_1:= \underline{r} + \varepsilon \le q+\varepsilon$. We also choose $\ell$ sufficiently large and $(r_j)$ sufficiently close to each other so that $r_j \, \varrho(r_{j-1}) < q+\varepsilon$ for all $2\le j\le \ell+1$. Thus, $\P$-almost surely for all large $i$,
$$
\sum_{x\in \T_i} \pi(x) \le (r_1)^i \deg^i + \sum_{j=2}^{\ell+1} i^2 \, \deg^i (q+\varepsilon)^i = (r_1)^i \deg^i + \ell \, i^2 \, \deg^i (q+\varepsilon)^i,
$$

\noindent which implies (recall: $\deg(q+\varepsilon)<1$) that $\sum_{i=0}^{n-1} {1\over \sum_{x\in \T_i} \pi(x)} \ge {c_2\over n^2\, \deg^n (q+\varepsilon)^n}$. Plugging this into (\ref{liggett}) yields that, $\P$-almost surely for all large $n$,
$$
P_\omega (\tau_n< \tau_0) \le c_3(\omega)\, n^2\, \deg^n (q+\varepsilon)^n \le [(q+2\varepsilon)\deg]^n.
$$

\noindent In particular, by writing $L(\tau_n):= \# \{ 1\le i \le \tau_n: \, X_i = e\}$, we obtain:
$$
P_\omega \left\{ L(\tau_n) \ge j \right\} = \left[ P_\omega (\tau_n> \tau_0) \right]^j \ge \left\{ 1- [(q+2\varepsilon)\deg]^n \right\}^j ,
$$

\noindent which, by the Borel--Cantelli lemma, yields that, $\P$-almost surely for all large $n$,
$$
L(\tau_n) \ge {1\over [(q+3\varepsilon) \deg]^n} , \qquad \hbox{\rm $P_\omega$-a.s.}
$$

\noindent Since $\{ L(\tau_n) \ge j \} \subset \{ \max_{0\le k \le 2j} |X_k| < n\}$, and since $\varepsilon$ can be as close to 0 as possible, we obtain the upper bound in Theorem \ref{t:posrec}.\hfill$\Box$

\bigskip

\noindent {\it Proof of Theorem \ref{t:posrec}: lower bound.} Assume $p< {1\over \deg}$. Recall that in this case, we have $\overline{r}<1$. Let $\varepsilon>0$ be small. Let $r \in (\underline{r}, \, \overline{r})$ be such that $\varrho(r) > {1\over \deg} \ee^\varepsilon$ and that $r\varrho(r) \ge q\ee^{-\varepsilon}$. Let $L$ be a large integer with $\deg^{-1/L} \ge \ee^{-\varepsilon}$ and satisfying (\ref{GW}) below.

We start by constructing a Galton--Watson tree $\G$, which is a certain subtree of $\T$. The first generation of $\G$, denoted by $\G_1$ and defined below, consists of vertices $x\in \T_L$ satisfying a certain property. The second generation of $\G$ is formed by applying the same procedure to each element of $\G_1$, and so on. To be precise,
$$
\G_1 = \G_1 (L,r) := \left\{ x\in \T_L: \, \min_{z\in ]\! ] e, \, x ]\! ]} \prod_{y\in ]\! ] e, \, z]\! ]} A(y) \ge r^L \right\} ,
$$

\noindent where $]\! ]e, \, x ]\! ]$ denotes as before the set of vertices (excluding $e$) lying on the shortest path relating $e$ and $x$. More generally, if $\G_i$ denotes the $i$-th generation of $\G$, then
$$
\G_{n+1} := \bigcup_{u\in \G_n } \left\{ x\in \T_{(n+1)L}: \, \min_{z\in ]\! ] u, \, x ]\! ]} \prod_{y\in ]\! ] u, \, z]\! ]} A(y) \ge r^L \right\} , \qquad n=1,2, \dots
$$

We claim that it is possible to choose $L$ sufficiently large such that
\begin{equation}
    \E(\# \G_1) \ge \ee^{-\varepsilon L} \deg^L
    \varrho(r)^L .
    \label{GW}
\end{equation}

\noindent Note that $\ee^{-\varepsilon L} \deg^L \varrho(r)^L>1$, since $\varrho(r) > {1\over \deg} \ee^\varepsilon$.

We admit (\ref{GW}) for the moment, which implies that $\G$ is super-critical. By theory of branching processes (Harris~\cite{harris}, p.~13), when $n$ goes to infinity, ${\# \G_{n/L} \over [\E(\# \G_1)]^{n/L} }$ converges almost surely (and in $L^2$) to a limit $W$ with $\P(W>0)>0$. Therefore, on the event $\{ W>0\}$, for all large $n$,
\begin{equation}
    \# (\G_{n/L}) \ge c_4(\omega)
    [\E(\# \G_1)]^{n/L}.
    \label{GnL}
\end{equation}

\noindent (For notational simplification, we only write our argument for the case when $n$ is a multiple of $L$. It is clear that our final conclusion holds for all large $n$.)

Recall that according to the Dirichlet principle (Griffeath and Liggett~\cite{griffeath-liggett}),
\begin{eqnarray}
    2\pi(e) P_\omega \left\{ \tau_n < \tau_0
    \right\}
 &=&\inf_{h: \, h(e)=1, \, h(z)=0, \, \forall |z|
    \ge n} \sum_{x,y\in \T} \pi(x) \omega(x,y)
    (h(x)- h(y))^2
    \nonumber
    \\
 &\ge& c_5\, \inf_{h: \, h(e)=1, \, h(z)=0, \,
    \forall z\in \T_n} \sum_{|x|<n} \sum_{y: \, x=
    {\buildrel \leftarrow \over y}} \pi(x) (h(x)-
    h(y))^2,
    \label{durrett}
\end{eqnarray}

\noindent the last inequality following from ellipticity condition on the environment. Clearly,
\begin{eqnarray*}
    \sum_{|x|<n} \sum_{y: \, x= {\buildrel
    \leftarrow \over y}} \pi(x) (h(x)- h(y))^2
 &=&\sum_{i=0}^{(n/L)-1} \sum_{x: \, iL \le |x| <
    (i+1) L} \sum_{y: \, x= {\buildrel \leftarrow
    \over y}} \pi(x) (h(x)- h(y))^2
    \\
 &:=&\sum_{i=0}^{(n/L)-1} I_i,
\end{eqnarray*}

\noindent with obvious notation. For any $i$,
$$
I_i \ge \deg^{-L} \sum_{v\in \G_{i+1}} \, \sum_{x\in [\! [ v^\uparrow, v[\! [} \, \sum_{y: \, x= {\buildrel \leftarrow \over y}} \pi(x) (h(x)- h(y))^2,
$$

\noindent where $v^\uparrow \in \G_i$ denotes the unique element of $\G_i$ lying on the path $[ \! [ e, v ]\! ]$ (in words, $v^\uparrow$ is the parent of $v$ in the Galton--Watson tree $\G$), and the factor $\deg^{-L}$ comes from the fact that each term $\pi(x) (h(x)- h(y))^2$ is counted at most $\deg^L$ times in the sum on the right-hand side.

By (\ref{pi}), for $x\in [\! [ v^\uparrow, v[\! [$, $\pi(x) \ge c_0 \, \prod_{u\in ]\! ]e, x]\! ]} A(u)$, which, by the definition of $\G$, is at least $c_0 \, r^{(i+1)L}$. Therefore,
\begin{eqnarray*}
    I_i
 &\ge& c_0 \, \deg^{-L} \sum_{v\in \G_{i+1}} \, \sum_{x\in [\!
    [ v^\uparrow, v[\! [} \, \sum_{y: \, x=
    {\buildrel \leftarrow \over y}} r^{(i+1)L}
    (h(x)- h(y))^2
    \\
 &\ge&c_0 \, \deg^{-L} r^{(i+1)L} \sum_{v\in \G_{i+1}} \,
    \sum_{y\in ]\! ] v^\uparrow, v]\! ]}
    (h({\buildrel \leftarrow \over y})- h(y))^2 .
\end{eqnarray*}

\noindent By the Cauchy--Schwarz inequality, $\sum_{y\in ]\! ] v^\uparrow, v]\! ]} (h({\buildrel \leftarrow \over y})- h(y))^2 \ge {1\over L} (h(v^\uparrow)-h(v))^2$. Accordingly,
$$
I_i \ge c_0 \, {\deg^{-L} r^{(i+1)L}\over L} \sum_{v\in \G_{i+1}} (h(v^\uparrow)-h(v))^2 ,
$$

\noindent which yields
\begin{eqnarray*}
    \sum_{i=0}^{(n/L)-1} I_i
 &\ge& c_0 \, {\deg^{-L}\over L} \sum_{i=0}^{(n/L)-1}
    r^{(i+1)L} \sum_{v\in \G_{i+1}} (h(v^\uparrow)-
    h(v))^2
    \\
 &\ge& c_0 \, {\deg^{-L}\over L} \deg^{-n/L} \sum_{v\in \G_{n/L}}
    \sum_{i=0}^{(n/L)-1} r^{(i+1)L}  (h(v^{(i)})-
    h(v^{(i+1)}))^2 ,
\end{eqnarray*}

\noindent where, $e=: v^{(0)}$, $v^{(1)}$, $v^{(2)}$, $\cdots$, $v^{(n/L)} := v$, is the shortest path (in $\G$) from $e$ to $v$, and the factor $\deg^{-n/L}$ results from the fact that each term $r^{(i+1)L}  (h(v^{(i)})- h(v^{(i+1)}))^2$ is counted at most $\deg^{n/L}$ times in the sum on the right-hand side.

By the Cauchy--Schwarz inequality, for all $h: \T\to \r$ with $h(e)=1$ and $h(z)=0$ ($\forall z\in \T_n$), we have
\begin{eqnarray*}
    \sum_{i=0}^{(n/L)-1} r^{(i+1)L}  (h(v^{(i)})-
    h(v^{(i+1)}))^2
 &\ge&{1\over \sum_{i=0}^{(n/L)-1} r^{-(i+1)L}} \,
    \left( \sum_{i=0}^{(n/L)-1} (h(v^{(i)})-
    h(v^{(i+1)})) \right)^{\! \! 2}
    \\
 &=&{1\over \sum_{i=0}^{(n/L)-1} r^{-(i+1)L}} \ge
    c_6 \, r^n.
\end{eqnarray*}

\noindent Therefore,
$$
\sum_{i=0}^{(n/L)-1} I_i \ge c_0c_6 \, r^n \, {\deg^{-L}\over L} \deg^{-n/L} \# (\G_{n/L}) \ge c_0 c_6 c_4(\omega) \, r^n \, {\deg^{-L}\over L} \deg^{-n/L} \, [\E (\# \G_1)]^{n/L}\, {\bf 1}_{ \{ W>0 \} },
$$

\noindent the last inequality following from (\ref{GnL}). Plugging this into (\ref{durrett}) yields that for all large $n$,
$$
P_\omega \left\{ \tau_n < \tau_0 \right\} \ge c_7(\omega) \, r^n \, {\deg^{-L}\over L} \deg^{-n/L} \, [\E (\# \G_1)]^{n/L}\, {\bf 1}_{ \{ W>0 \} } .
$$

\noindent Recall from (\ref{GW}) that $\E(\# \G_1) \ge \ee^{-\varepsilon L} \deg^L \varrho(r)^L$. Therefore, on $\{W>0\}$, for all large $n$, $P_\omega \{ \tau_n < \tau_0 \} \ge c_8(\omega) (\ee^{-\varepsilon} \deg^{-1/L} \deg r \varrho(r))^n$, which is no smaller than $c_8(\omega) (\ee^{-3\varepsilon} q \deg)^n$ (since $\deg^{-1/L} \ge \ee^{-\varepsilon}$ and $r \varrho(r) \ge q \ee^{-\varepsilon}$ by assumption).
Thus, by writing $L(\tau_n) := \#\{ 1\le i\le n: \; X_i = e \}$ as before, we have, on $\{ W>0 \}$,
$$
P_\omega \left\{ L(\tau_n) \ge j \right\} = \left[ P_\omega (\tau_n> \tau_0) \right]^j \le [1- c_8(\omega) (\ee^{-3\varepsilon} q \deg)^n ]^j.
$$

\noindent By the Borel--Cantelli lemma, for $\P$-almost all $\omega$, on $\{W>0\}$, we have, $P_\omega$-almost surely for all large $n$, $L(\tau_n) \le 1/(\ee^{-4\varepsilon} q \deg)^n$, i.e.,
$$
\max_{0\le k\le \tau_0(\lfloor 1/(\ee^{-4\varepsilon} q \deg)^n\rfloor )} |X_k| \ge n ,
$$

\noindent where $0<\tau_0(1)<\tau_0(2)<\cdots$ are the successive return times to the root $e$ by the walk (thus $\tau_0(1) = \tau_0$). Since the walk is positive recurrent, $\tau_0(\lfloor 1/(\ee^{-4\varepsilon} q \deg)^n\rfloor ) \sim {1\over (\ee^{-4\varepsilon} q \deg)^n} E_\omega [\tau_0]$  (for $n\to \infty$), $P_\omega$-almost surely ($a_n \sim b_n$ meaning $\lim_{n\to \infty}Ê{a_n \over b_n} =1$). Therefore, for $\P$-almost all $\omega \in \{ W>0\}$,
$$
\liminf_{n\to \infty} {\max_{0\le k\le n} |X_k| \over \log n} \ge {1\over \log[1/(q\deg)]}, \qquad \hbox{\rm $P_\omega$-a.s.}
$$

\noindent Recall that $\P\{ W>0\}>0$. Since modifying a finite number of transition probabilities does not change the value of $\liminf_{n\to \infty} {\max_{0\le k\le n} |X_k| \over \log n}$, we obtain the lower bound in Theorem \ref{t:posrec}.

It remains to prove (\ref{GW}). Let $(A^{(i)})_{i\ge 1}$ be an i.i.d.\ sequence of random variables distributed as $A$. Clearly, for any $\delta\in (0,1)$,
\begin{eqnarray*}
    \E( \# \G_1)
 &=& \deg^L \, \P\left( \, \sum_{i=1}^\ell \log
    A^{(i)} \ge L \log r , \, \forall 1\le \ell \le
    L\right)
    \\
 &\ge& \deg^L \, \P \left( \, (1-\delta) L
    \log r \ge \sum_{i=1}^\ell \log A^{(i)} \ge L
    \log r , \, \forall 1\le \ell \le L\right) .
\end{eqnarray*}

\noindent We define a new probability $\Q$ by
$$
{\mathrm{d} \Q \over \mathrm{d}\P} := {\ee^{t \log A} \over \E(\ee^{t \log A})} = {A^t \over \E(A^t)},
$$

\noindent for some $t\ge 0$. Then
\begin{eqnarray*}
    \E(\# \G_1)
 &\ge& \deg^L \, \E_\Q \left[ \, {[\E(A^t)]^L \over
    \exp\{ t \sum_{i=1}^L \log A^{(i)}\} }\,
    {\bf 1}_{\{ (1-\delta) L \log r \ge
    \sum_{i=1}^\ell \log A^{(i)} \ge L \log r , \,
    \forall 1\le \ell \le L\} } \right]
    \\
 &\ge& \deg^L \, {[\E(A^t)]^L \over r^{t (1-
    \delta) L} } \, \Q \left( (1-
    \delta) L \log r \ge \sum_{i=1}^\ell \log
    A^{(i)} \ge L \log r , \, \forall 1\le \ell \le
    L \right).
\end{eqnarray*}

\noindent To choose an optimal value of $t$, we fix $\widetilde{r}\in (r, \, \overline{r})$ with $\widetilde{r} < r^{1-\delta}$. Our choice of $t=t^*$ is such that $\varrho(\widetilde{r}) = \inf_{t\ge 0} \{ \widetilde{r}^{-t} \E(A^t)\} = \widetilde{r}^{-t^*} \E(A^{t^*})$. With this choice, we have $\E_\Q(\log A)=\log \widetilde{r}$, so that $\Q \{ (1- \delta) L \log r \ge \sum_{i=1}^\ell \log A^{(i)} \ge L \log r , \, \forall 1\le \ell \le L\} \ge c_9$. Consequently,
$$
\E(\# \G_1) \ge c_9 \, \deg^L \, {[\E(A^{t^*})]^L \over r^{t^* (1- \delta) L} }= c_9 \, \deg^L \, {[ \widetilde{r}^{\,t^*} \varrho(\widetilde{r})]^L \over r^{t^* (1- \delta) L} } \ge c_9 \, r^{\delta t^* L}  \deg^L \varrho(\widetilde{r})^L .
$$

\noindent Since $\delta>0$ can be as close to $0$ as possible, the continuity of $\varrho(\cdot)$ on $[\underline{r}, \, \overline{r})$ yields (\ref{GW}), and thus completes the proof of Theorem \ref{t:posrec}.\hfill$\Box$

\section{Some elementary inequalities}
\label{s:proba}

We collect some elementary inequalities in this section. They will be of use in the next sections, in the study of the null recurrence case.

\medskip

\begin{lemma}
\label{l:exp}
Let $\xi\ge 0$ be a random variable.

 {\rm (i)} Assume that $\e(\xi^a)<\infty$ for some
 $a>1$. Then for any $x\ge 0$,
 \begin{equation}
    {\e[({\xi\over x+\xi})^a] \over [\e ( {\xi\over
    x+\xi})]^a} \le {\e (\xi^a) \over [\e \xi]^a} .
    \label{RSD}
 \end{equation}

 {\rm (ii)} If $\e (\xi) < \infty$, then for any $0
 \le \lambda \le 1$ and $t \ge 0$,
 \begin{equation}
    \e \left\{ \exp \left( - t\, { (\lambda+\xi)/
    (1+\xi) \over \e [(\lambda+\xi)/
    (1+\xi)] } \right) \right\} \le \e \left\{
    \exp\left( - t\, { \xi \over \e (\xi)}
    \right) \right\} .
    \label{exp}
 \end{equation}

\end{lemma}

\medskip

\noindent {\bf Remark.} When $a=2$, (\ref{RSD}) is a special case of Lemma 6.4 of Pemantle and Peres~\cite{pemantle-peres2}.

\bigskip

\noindent {\it Proof of Lemma \ref{l:exp}.} We actually prove a very general result, stated as follows. Let $\varphi : (0, \infty) \to \r$ be a convex ${\cal C}^1$-function. Let $x_0 \in \r$ and let $I$ be an open interval containing $x_0$. Assume that $\xi$ takes values in a Borel set $J \subset \r$ (for the moment, we do not assume $\xi\ge 0$). Let $h: I \times J \to (0, \infty)$ and ${\partial h\over \partial x}: I \times J \to \r$ be measurable functions such that
\begin{itemize}

 \item $\e \{ h(x_0, \xi)\} <\infty$ and $\e \{
       |\varphi ({ h(x_0,\xi) \over \e h(x_0, \xi)}
       )| \} < \infty$;

 \item $\e[\sup_{x\in I} \{ | {\partial h\over
       \partial x} (x, \xi)| + |\varphi' ({h(x,
       \xi) \over \e h(x, \xi)} ) | \,
       ({| {\partial h\over \partial x} (x, \xi) |
       \over \e \{ h(x, \xi)\} } + {h(x, \xi) \over
       [\e \{ h(x, \xi)\}]^2 } | \e \{ {\partial
       h\over \partial x} (x, \xi) \} | )\} ] <
       \infty$;

 \item both $y \to h(x_0, y)$ and $y \to {
       \partial \over \partial x} \log
       h(x,y)|_{x=x_0}$ are monotone on $J$.

\end{itemize}

\noindent Then
\begin{equation}
    {\d \over \d x} \e \left\{ \varphi\left({
    h(x,\xi) \over \e h(x, \xi)}\right) \right\}
    \Big|_{x=x_0} \ge 0, \qquad \hbox{\rm or}\qquad
    \le 0,
    \label{monotonie}
\end{equation}

\noindent depending on whether $h(x_0, \cdot)$ and ${\partial \over \partial x} \log h(x_0,\cdot)$ have the same monotonicity.

To prove (\ref{monotonie}), we observe that by the integrability assumptions,
\begin{eqnarray*}
 & &{\d \over \d x} \e \left\{ \varphi\left({
    h(x,\xi) \over \e h(x,\xi)}\right) \right\}
    \Big|_{x=x_0}
    \\
 &=&{1 \over ( \e h(x_0, \xi))^2}\, \e \left(
    \varphi'( h(x_0, \xi) ) \left[ {\partial h \over
    \partial x} (x_0, \xi) \e h(x_0, \xi) -
    h(x_0, \xi) \e {\partial h \over \partial x}
    (x_0, \xi) \right] \right) .
\end{eqnarray*}

\noindent Let $\widetilde \xi$ be an independent copy of $\xi$. The expectation expression $\e(\varphi'( h(x_0, \xi) ) [\cdots])$ on the right-hand side is
\begin{eqnarray*}
 &=& \e \left(
    \varphi'( h(x_0, \xi) ) \left[ {\partial h \over
    \partial x} (x_0, \xi)   h(x_0, \widetilde\xi) -
    h(x_0, \xi)   {\partial h \over \partial x}
    (x_0, \widetilde\xi) \right] \right)
    \\
 &=& {1 \over 2}\, \e \left(
    \left[ \varphi'( h(x_0, \xi) ) - \varphi'(
    h(x_0, \widetilde\xi) )\right]
    \left[ {\partial h \over \partial x} (x_0, \xi)
    h(x_0, \widetilde\xi) - h(x_0, \xi) {\partial
    h \over \partial x} (x_0, \widetilde\xi)
    \right] \right)
    \\
 &=& {1 \over 2}\, \e \left( h(x_0, \xi) h(x_0,
    \widetilde \xi) \, \eta \right) ,
\end{eqnarray*}

\noindent where
$$
\eta := \left[ \varphi'( h(x_0, \xi) ) - \varphi'(  h(x_0, \widetilde\xi) ) \right] \, \left[ {\partial  \log h \over \partial x} (x_0, \xi) - {\partial \log h \over \partial x} (x_0, \widetilde\xi) \right] .
$$

\noindent Therefore,
$$
{\d \over \d x} \e \left\{ \varphi\left({ h(x,\xi) \over \e h(x,\xi)}\right) \right\} \Big|_{x=x_0}
\; = \; {1 \over 2( \e h(x_0, \xi))^2}\, \e \left( h(x_0, \xi) h(x_0, \widetilde \xi) \, \eta \right) .
$$

\noindent Since $\eta \ge 0$ or $\le 0$ depending on whether $h(x_0, \cdot)$ and ${\partial \over \partial x} \log h(x_0,\cdot)$ have the same monotonicity, this yields (\ref{monotonie}).

To prove (\ref{RSD}) in Lemma \ref{l:exp}, we take $x_0\in (0,\, \infty)$, $J= \r_+$, $I$ a finite open interval containing $x_0$ and away from 0, $\varphi(z)= z^a$, and $h(x,y)= { y \over x+ y}$, to see that the function $x\mapsto {\e[({\xi\over x+\xi})^a] \over [\e ( {\xi\over x+\xi})]^a}$ is non-decreasing on $(0, \infty)$.  By dominated convergence,
$$
\lim_{x \to\infty}  {\e[({\xi\over x+\xi})^a] \over [\e ( {\xi\over x+\xi})]^a}=  \lim_{x \to\infty}  {\e[({\xi\over 1+\xi/x})^a] \over [\e ( {\xi\over 1+\xi/x})]^a} =  {\e (\xi^a) \over [\e \xi]^a} ,
$$

\noindent yielding (\ref{RSD}).

The proof of (\ref{exp}) is similar. Indeed,  applying (\ref{monotonie}) to the functions $\varphi(z)= \ee^{-t z}$ and $ h(x, y) = {x + y \over 1+ y}$ with $x\in (0,1)$, we get that
the function $x \mapsto \e \{ \exp ( - t {
(x+\xi)/(1+\xi) \over \e [(x+\xi)/(1+\xi)]} )\}$ is non-increasing on $(0,1)$; hence for $\lambda \in [0,\, 1]$,
$$
\e \left\{ \exp \left( - t { (\lambda+\xi)/(1+\xi) \over \e [(\lambda+\xi)/(1+\xi)] } \right) \right\} \le \e \left\{ \exp \left( - t { \xi /(1+\xi) \over \e [\xi/(1+\xi)] } \right) \right\}.
$$

\noindent On the other hand, we take $\varphi(z)= \ee^{-t z}$ and $h(x,y) = {y \over 1+ xy}$ (for $x\in (0, 1)$) in (\ref{monotonie}) to see that $x \mapsto \e \{ \exp ( - t { \xi /(1+x \xi) \over \e [\xi /(1+x\xi)] } ) \}$ is non-increasing on $(0,1)$. Therefore,
$$
\e \left\{ \exp \left( - t { \xi /(1+\xi) \over \e [\xi/(1+\xi)] } \right) \right\} \le \e \left\{ \exp\left( - t \, { \xi \over \e (\xi)}\right) \right\} ,
$$

\noindent which implies (\ref{exp}).\hfill$\Box$

\medskip

\begin{lemma}
 \label{l:moment}
 Let $\xi_1$, $\cdots$, $\xi_k$ be independent
 non-negative random variables such that for some
 $a\in [1,\, 2]$, $\e(\xi_i^a)<\infty$ $(1\le i\le
 k)$. Then
 $$
 \e \left[ (\xi_1 + \cdots + \xi_k)^a \right] \le
 \sum_{k=1}^k \e(\xi_i^a) + (k-1) \left(
 \sum_{i=1}^k \e \xi_i \right)^a.
 $$

\end{lemma}

\medskip

\noindent {\it Proof.} By induction on $k$, we only need to prove the lemma in case $k=2$. Let
$$
h(t) := \e \left[ (\xi_1 + t\xi_2)^a \right] - \e(\xi_1^a) - t^a \e(\xi_2^a) - (\e \xi_1 + t \e \xi_2)^a, \qquad t\in [0,1].
$$

\noindent Clearly, $h(0) = - (\e \xi_1)^a \le 0$. Moreover,
$$
h'(t) = a \e \left[ (\xi_1 + t\xi_2)^{a-1} \xi_2 \right] - a t^{a-1} \e(\xi_2^a)  - a(\e \xi_1 + t \e \xi_2)^{a-1} \e(\xi_2) .
$$

\noindent Since $(x+y)^{a-1} \le x^{a-1} + y^{a-1}$ (for $1\le a\le 2$), we have
\begin{eqnarray*}
    h'(t)
 &\le& a \e \left[ (\xi_1^{a-1} + t^{a-1}\xi_2^{a
    -1}) \xi_2 \right] - a t^{a-1} \e(\xi_2^a)  -
    a(\e \xi_1)^{a-1} \e(\xi_2)
    \\
 &=& a \e (\xi_1^{a-1}) \e(\xi_2) - a(\e \xi_1)^{a
    -1} \e(\xi_2) \le 0,
\end{eqnarray*}

\noindent by Jensen's inequality (for $1\le a\le 2$). Therefore, $h \le 0$ on $[0,1]$. In particular, $h(1) \le 0$, which implies Lemma \ref{l:moment}.\hfill$\Box$

\bigskip

The following inequality, borrowed from page 82 of Petrov~\cite{petrov}, will be of frequent use.

\medskip

\begin{fact}
 \label{f:petrov}
 Let $\xi_1$, $\cdots$, $\xi_k$ be independent
 random variables. We assume that for any $i$,
 $\e(\xi_i)=0$ and $\e(|\xi_i|^a) <\infty$, where
 $1\le a\le 2$. Then
 $$
 \e \left( \, \left| \sum_{i=1}^k \xi_i \right| ^a
 \, \right) \le 2 \sum_{i=1}^k \e( |\xi_i|^a).
 $$

\end{fact}

\medskip

\begin{lemma}
 \label{l:abc}
 Fix $a >1$. Let $(u_j)_{j\ge 1}$ be a sequence of
 positive numbers, and let $(\lambda_j)_{j\ge 1}$ be
 a sequence of non-negative numbers.

 {\rm (i)} If there exists some constant $c_{10}>0$
 such that for all $n\ge 2$,
 $$
 u_{j+1} \le \lambda_n + u_j - c_{10}\, u_j^{a},
 \qquad \forall 1\le j \le n-1,
 $$
 then we can find a constant $c_{11}>0$ independent
 of $n$ and $(\lambda_j)_{j\ge 1}$, such that
 $$
 u_n \le c_{11} \, ( \lambda_n^{1/a} +
 n^{- 1/(a-1)}), \qquad \forall n\ge 1.
 $$

 {\rm (ii)} Fix $K>0$. Assume that
 $\lim_{j\to\infty} u_j=0$ and that $\lambda_n \in
 [0, \, {K\over n}]$ for all $n\ge 1$.
 If there exist $c_{12}>0$ and $c_{13}>0$ such that
 for all $n\ge 2$,
 $$
 u_{j+1} \ge \lambda_n + (1- c_{12} \lambda_n) u_j -
 c_{13} \, u_j^a , \qquad \forall 1 \le j \le n-1,
 $$
 then for some $c_{14}>0$ independent of $n$ and
 $(\lambda_j)_{j\ge 1}$ $(c_{14}$ may depend on
 $K)$,
 $$
 u_n \ge c_{14} \, ( \lambda_n^{1/a} + n^{- 1/(a-1)}
 ), \qquad \forall n\ge 1.
 $$

\end{lemma}

\medskip

{\noindent\it Proof.} (i) Put $\ell = \ell(n) := \min\{n, \, \lambda_n^{- (a-1)/a} \}$. There are two possible situations.

First situation: there exists some $j_0 \in [n- \ell,  n-1]$ such that $u_{j_0} \le ({2 \over c_{10}})^{1/a}\,   \lambda_n^{1/a}$. Since $u_{j+1} \le \lambda_n + u_j$ for all $j\in [j_0, n-1]$, we have
$$
u_n \le (n-j_0 ) \lambda_n + u_{j_0} \le  \ell \lambda_n + ({2 \over c_{10}})^{1/a}\, \lambda_n^{1/a} \le (1+ ({2 \over c_{10}})^{1/a})\, \lambda_n^{1/a},
$$

\noindent which implies the desired upper bound.

Second situation: $u_j > ({2 \over c_{10}})^{1/a}\, \lambda_n^{1/a}$, $\forall \, j \in [n- \ell,  n-1]$. Then $c_{10}\, u_j^{a} > 2\lambda_n$, which yields
$$
u_{j+1} \le u_j - {c_{10} \over 2} u_j^a, \qquad \forall \, j \in [n- \ell,  n-1].
$$

\noindent Since $a>1$ and $(1-y)^{1-a} \ge 1+ (a-1) y$ (for $0< y< 1$), this yields, for $j \in [n- \ell,  n-1]$,
$$
u_{j+1}^{1-a} \ge u_j^{1-a} \, \left( 1 - {c_{10} \over 2} u_j^{a-1} \right)^{ 1-a}  \ge  u_j^{ 1-a} \, \left( 1 + {c_{10} \over 2} (a-1)\, u_j^{a-1} \right) = u_j^{1-a} + {c_{10} \over 2} (a-1) .
$$

\noindent Therefore, $u_n^{1-a} \ge c_{15}\, \ell$ with $c_{15}:= {c_{10} \over 2} (a-1)$. As a consequence, $u_n \le (c_{15}\, \ell)^{- 1/(a-1)} \le (c_{15})^{- 1/(a-1)} \, ( n^{- 1/(a-1)} + \lambda_n^{1/a} )$, as desired.

(ii) Let us first prove:
\begin{equation}
    \label{c7}
    u_n \ge c_{16}\, n^{- 1/(a-1)}.
\end{equation}

To this end,  let $n$ be large and define $v_j := u_j \, (1- c_{12} \lambda_n)^{ -j} $ for $1 \le j \le n$. Since $u_{j+1} \ge (1- c_{12} \lambda_n) u_j - c_{13} u_j^a $ and $\lambda_n \le K/n$, we get
$$
v_{j+1} \ge v_j - c_{13} (1- c_{12} \lambda_n)^{(a-1)j-1}\, v_j^a\ge  v_j - c_{17} \, v_j^a, \qquad \forall\, 1\le j \le n-1.
$$

\noindent Since $u_j \to 0$, there exists some $j_0>0$ such that for all $n>j \ge j_0$, we have $c_{17} \, v_j^{a-1} < 1/2$, and
$$
v_{j+1}^{1-a} \le v_j^{1-a}\, \left( 1- c_{17} \, v_j^{a-1}\right)^{1-a} \le  v_j^{1-a}\, \left( 1+ c_{18} \, v_j^{a-1}\right) = v_j^{1-a} + c_{18}.
$$

\noindent It follows that $v_n^{1-a} \le c_{18}\, (n-j_0) + v_{j_0}^{1-a}$, which implies (\ref{c7}).

It remains to show that $u_n \ge c_{19} \, \lambda_n^{1/a}$. Consider a large $n$. The function $h(x):= \lambda_n + (1- c_{12} \lambda_n) x - c_{13} x^a$  is increasing on $[0, c_{20}]$ for some fixed constant $c_{20}>0$. Since $u_j \to 0$, there exists
$j_0$ such that $u_j \le c_{20}$ for all $j \ge j_0$. We claim there exists $j \in [j_0, n-1]$ such that $u_j > ({\lambda_n\over 2c_{13}})^{1/a}$: otherwise, we would have $c_{13}\, u_j^a \le {\lambda_n\over 2} \le \lambda_n$ for all $j \in [j_0, n-1]$, and thus
$$
u_{j+1} \ge (1- c_{12}\, \lambda_n) u_j \ge \cdots \ge (1- c_{12}\,\lambda_n)^{j-j_0} \, u_{j_0} ;
$$

\noindent in particular, $u_n \ge (1- c_{12}\, \lambda_n)^{n-j_0} \, u_{j_0}$  which would contradict the assumption $u_n \to 0$ (since $\lambda_n \le K/n$).

Therefore, $u_j > ({\lambda_n\over 2c_{13}})^{1/a}$ for some $j\ge j_0$. By monotonicity of $h(\cdot)$ on $[0, c_{20}]$,
$$
u_{j+1} \ge h(u_j) \ge h\left(({\lambda_n\over 2
c_{13}})^{1/a}\right) \ge ({\lambda_n\over 2 c_{13}})^{1/a},
$$

\noindent the last inequality being elementary. This leads to: $u_{j+2} \ge h(u_{j+1}) \ge h(({\lambda_n\over 2 c_{13}})^{1/a} ) \ge ({\lambda_n\over 2 c_{13}})^{1/a}$. Iterating the procedure, we obtain: $u_n \ge ({\lambda_n\over 2 c_{13}})^{1/a}$ for all $n> j_0$, which completes the proof of the Lemma.\hfill$\Box$

\section{Proof of Theorem \ref{t:nullrec}}
\label{s:nullrec}

Let $n\ge 2$, and let as before
$$
\tau_n := \inf\left\{ i\ge 1: X_i \in \T_n
\right\} .
$$

\noindent We start with a characterization of the distribution of $\tau_n$ via its Laplace transform $\e ( \ee^{- \lambda \tau_n} )$, for $\lambda \ge 0$. To state the result, we define $\alpha_{n,\lambda}(\cdot)$, $\beta_{n,\lambda}(\cdot)$ and $\gamma_n(\cdot)$ by $\alpha_{n,\lambda}(x) = \beta_{n,\lambda} (x) = 1$ and $\gamma_n(x)=0$ (for $x\in \T_n$), and
\begin{eqnarray}
    \alpha_{n,\lambda}(x)
 &=& \ee^{-\lambda} \, {\sum_{i=1}^\deg A(x_i)
    \alpha_{n,\lambda} (x_i) \over 1+
    \sum_{i=1}^\deg A(x_i) \beta_{n,\lambda}
    (x_i)},
    \label{alpha}
    \\
    \beta_{n,\lambda}(x)
 &=& {(1-\ee^{-2\lambda}) + \sum_{i=1}^\deg A(x_i)
    \beta_{n,\lambda} (x_i) \over 1+
    \sum_{i=1}^\deg A(x_i) \beta_{n,\lambda} (x_i)},
    \label{beta}
    \\
    \gamma_n(x)
 &=& {[1/\omega(x, {\buildrel \leftarrow \over x} )]
    + \sum_{i=1}^\deg A(x_i) \gamma_n(x_i) \over 1+
    \sum_{i=1}^\deg A(x_i) \beta_n(x_i)} , \qquad
    1\le |x| < n,
    \label{gamma}
\end{eqnarray}

\noindent where  $\beta_n(\cdot) := \beta_{n,0}(\cdot)$, and for any $x\in \T$, $\{x_i\}_{1\le i\le \deg}$ stands as before for the set of children of $x$.

\medskip

\begin{proposition}
 \label{p:tau}
 We have, for $n\ge 2$,
 \begin{eqnarray}
     E_\omega\left( \ee^{- \lambda \tau_n} \right)
  &=&\ee^{-\lambda} \, {\sum_{i=1}^\deg \omega (e,
     e_i) \alpha_{n,\lambda} (e_i) \over
     \sum_{i=1}^\deg \omega (e, e_i)
     \beta_{n,\lambda} (e_i)}, \qquad \forall
     \lambda \ge 0,
     \label{Laplace-tau}
     \\
     E_\omega(\tau_n)
  &=& {1+ \sum_{i=1}^\deg \omega(e,e_i) \gamma_n
     (e_i) \over \sum_{i=1}^\deg \omega(e,e_i)
     \beta_n(e_i)}.
     \label{E(tau)}
 \end{eqnarray}

\end{proposition}

\medskip

\noindent {\it Proof of Proposition \ref{p:tau}.} Identity (\ref{E(tau)}) can be found in Rozikov~\cite{rozikov}. The proof of (\ref{Laplace-tau}) is along similar lines; so we feel free to give an outline only. Let $g_{n, \lambda}(x) := E_\omega (\ee^{- \lambda \tau_n} \, | \, X_0=x)$. By the Markov property, $g_{n, \lambda}(x) = \ee^{-\lambda} \sum_{i=1}^\deg \omega(x, x_i)g_{n, \lambda}(x_i) + \ee^{-\lambda} \omega(x, {\buildrel \leftarrow \over x}) g_{n, \lambda}({\buildrel \leftarrow \over x})$, for $|x| < n$. By induction on $|x|$ (such that $1\le |x| \le n-1$), we obtain: $g_{n, \lambda}(x) = \ee^\lambda (1- \beta_{n, \lambda} (x)) g_{n, \lambda}({\buildrel \leftarrow \over x}) + \alpha_{n, \lambda} (x)$, from which (\ref{Laplace-tau}) follows.

Probabilistic interpretation: for $1\le |x| <n$, if $T_{\buildrel \leftarrow \over x} := \inf \{ k\ge 0: X_k= {\buildrel \leftarrow \over x} \}$, then $\alpha_{n, \lambda} (x) = E_\omega [ \ee^{-\lambda \tau_n} {\bf 1}_{ \{ \tau_n < T_{\buildrel \leftarrow \over x} \} } \, | \, X_0=x]$, $\beta_{n, \lambda} (x) = 1- E_\omega [ \ee^{-\lambda (1+ T_{\buildrel \leftarrow \over x}) } {\bf 1}_{ \{ \tau_n > T_{\buildrel \leftarrow \over x} \} } \, | \, X_0=x]$, and $\gamma_n (x) = E_\omega [ (\tau_n \wedge T_{\buildrel \leftarrow \over x}) \, | \, X_0=x]$. We do not use these identities in the paper.\hfill$\Box$

\bigskip

It turns out that $\beta_{n,\lambda}(\cdot)$ is closely related to Mandelbrot's multiplicative cascade~\cite{mandelbrot}. Let
\begin{equation}
    M_n := \sum_{x\in \T_n} \prod_{y\in ] \! ] e, \,
    x] \! ] } A(y) , \qquad n\ge 1,
    \label{Mn}
\end{equation}

\noindent where $] \! ] e, \,x] \! ]$ denotes as before the shortest path relating $e$ to $x$. We mention that $(A(e_i), \, 1\le i\le \deg)$ is a random vector independent of $(\omega(x,y), \, |x|\ge 1, \, y\in \T)$, and is distributed as $(A(x_i), \, 1\le i\le \deg)$, for any $x\in \T_m$ with $m\ge 1$.

Let us recall some properties of $(M_n)$ from Theorem 2.2 of Liu~\cite{liu00} and Theorem 2.5 of Liu~\cite{liu01}: under the conditions $p={1\over \deg}$ and $\psi'(1)<0$, $(M_n)$ is a martingale, bounded in $L^a$ for any $a\in [1, \kappa)$; in particular,
\begin{equation}
    M_\infty := \lim_{n\to \infty} M_n \in (0,
    \infty),
    \label{cvg-M}
\end{equation}

\noindent exists $\P$-almost surely and in $L^a(\P)$, and
\begin{equation}
    \E\left( \ee^{-s M_\infty} \right) \le
    \exp\left( - c_{21} \, s^{c_{22}}\right), \qquad
    \forall s\ge 1;
    \label{M-lowertail}
\end{equation}

\noindent furthermore, if $1<\kappa< \infty$, then we also have
\begin{equation}
    {c_{23}\over x^\kappa} \le \P\left( M_\infty >
    x\right) \le {c_{24}\over x^\kappa}, \qquad
    x\ge 1.
    \label{M-tail}
\end{equation}

We now summarize the asymptotic properties of $\beta_{n,\lambda}(\cdot)$ which will be needed later on.

\bigskip

\begin{proposition}
 \label{p:beta-gamma}
 Assume $p= {1\over \deg}$ and $\psi'(1)<0$.

 {\rm (i)} For any $1\le i\le \deg$, $n\ge 2$, $t\ge
 0$ and $\lambda \in [0, \, 1]$, we have
 \begin{equation}
     \E \left\{ \exp \left[ -t \, {\beta_{n,
     \lambda} (e_i) \over \E[\beta_{n, \lambda}
     (e_i)]} \right] \right\} \le \left\{\E \left(
     \ee^{-t\, M_n/\Theta} \right)
     \right\}^{1/\deg} ,
     \label{comp-Laplace}
 \end{equation}
 where, as before, $\Theta:= \hbox{\rm ess sup}(A) <
 \infty$.

 {\rm (ii)} If $\kappa\in (2, \infty]$, then for any
 $1\le i\le \deg$ and all $n\ge 2$ and $\lambda \in
 [0, \, {1\over n}]$,
 \begin{equation}
     c_{25} \left( \sqrt {\lambda} + {1\over n}
     \right) \le \E[\beta_{n, \lambda}(e_i)]
     \le c_{26} \left( \sqrt {\lambda} + {1\over
     n} \right).
     \label{E(beta):kappa>2}
 \end{equation}

 {\rm (iii)} If $\kappa\in (1,2]$, then for any
 $1\le i\le \deg$, when $n\to \infty$ and uniformly
 in $\lambda \in [0, {1\over n}]$,
 \begin{equation}
     \E[\beta_{n, \lambda}(e_i)] \; \approx \;
     \lambda^{1/\kappa} + {1\over n^{1/(\kappa-1)}}
     ,
     \label{E(beta):kappa<2}
 \end{equation}
 where $a_n \approx b_n$ denotes as before
 $\lim_{n\to \infty} \, {\log a_n \over \log b_n}
 =1$.

\end{proposition}

\bigskip

The proof of Proposition \ref{p:beta-gamma} is postponed until Section \ref{s:beta-gamma}. By admitting it for the moment, we are able to prove Theorem \ref{t:nullrec}.

\bigskip

\noindent {\it Proof of Theorem \ref{t:nullrec}.} Assume $p= {1\over \deg}$ and $\psi'(1)<0$.

Let $\pi$ be an invariant measure. By (\ref{pi}) and the definition of $(M_n)$, $\sum_{x\in \T_n} \pi(x) \ge c_0 \, M_n$. Therefore by (\ref{cvg-M}), we have $\sum_{x\in \T} \pi(x) =\infty$, $\P$-a.s., implying that $(X_n)$ is null recurrent.

We proceed to prove the lower bound in (\ref{nullrec}). By (\ref{gamma}) and the ellipticity condition on the environment, $\gamma_n (x) \le {1\over \omega(x, {\buildrel \leftarrow \over x} )} + \sum_{i=1}^\deg A(x_i) \gamma_n(x_i) \le c_{27} + \sum_{i=1}^\deg A(x_i) \gamma_n(x_i)$. Iterating the argument yields
$$
\gamma_n (e_i) \le c_{27} \left( 1+ \sum_{j=2}^{n-1} M_j^{(e_i)}\right), \qquad n\ge 3,
$$

\noindent where
$$
M_j^{(e_i)} := \sum_{x\in \T_j} \prod_{y\in ] \! ] e_i, x] \! ]} A(y).
$$

\noindent For future use, we also observe that
\begin{equation}
    \label{defMei1}
    M_n= \sum_{i=1}^\deg \, A(e_i) \, M^{(e_i)}_n,
    \qquad n\ge 2.
\end{equation}

Let $1\le i\le \deg$. Since $(M_j^{(e_i)}, \, j\ge 2)$ is distributed as $(M_{j-1}, \, j\ge 2)$, it follows from (\ref{cvg-M}) that $M_j^{(e_i)}$ converges (when $j\to \infty$) almost surely, which
implies $\gamma_n (e_i) \le c_{28}(\omega) \, n$.
Plugging this into (\ref{E(tau)}), we see that for all $n\ge 3$,
\begin{equation}
    E_\omega \left( \tau_n \right) \le {c_{29}(\omega) \, n \over \sum_{i=1}^\deg
    \omega(e,e_i) \beta_n(e_i)} \le {c_{30}(\omega)
    \, n \over \beta_n(e_1)},
    \label{toto2}
\end{equation}

\noindent the last inequality following from the ellipticity assumption on the environment.

We now bound $\beta_n(e_1)$ from below (for large $n$). Let $1\le i\le \deg$. By (\ref{comp-Laplace}), for $\lambda \in [0,\, 1]$ and $s\ge 0$,
$$
\E \left\{ \exp \left[ -s \, {\beta_{n, \lambda}
(e_i) \over \E [\beta_{n, \lambda} (e_i)]} \right] \right\} \le \left\{ \E \left( \ee^{-s \, M_n/\Theta} \right) \right\}^{1/\deg} \le \left\{ \E \left(\ee^{-s \, M_\infty/\Theta} \right) \right\}^{1/\deg} ,
$$

\noindent where, in the last inequality, we used the fact that $(M_n)$ is a uniformly integrable martingale. Let $\varepsilon>0$. Applying (\ref{M-lowertail}) to $s:= n^{\varepsilon}$, we see that
\begin{equation}
    \sum_n \E \left\{ \exp \left[ -n^{\varepsilon}
    {\beta_{n, \lambda} (e_i) \over \E[\beta_{n,
    \lambda} (e_i)]} \right] \right\} <\infty .
    \label{toto3}
\end{equation}

\noindent In particular, $\sum_n \exp [ -n^{\varepsilon} {\beta_n (e_1) \over \E [\beta_n (e_1)]} ]$ is $\P$-almost surely finite (by taking $\lambda=0$; recalling that $\beta_n (\cdot) := \beta_{n, 0} (\cdot)$). Thus, for $\P$-almost all $\omega$ and all sufficiently large $n$, $\beta_n (e_1) \ge n^{-\varepsilon} \, \E [\beta_n (e_1)]$. Going back to (\ref{toto2}), we see that for $\P$-almost all $\omega$ and all sufficiently large $n$,
$$
E_\omega \left( \tau_n \right) \le {c_{30}(\omega) \, n^{1+\varepsilon} \over \E [\beta_n (e_1)]}.
$$

\noindent Let $m(n):= \lfloor {n^{1+2\varepsilon} \over \E [\beta_n (e_1)]} \rfloor$. By Chebyshev's inequality, for $\P$-almost all $\omega$ and all sufficiently large $n$, $P_\omega ( \tau_n \ge m(n) ) \le c_{31}(\omega) \, n^{-\varepsilon}$. Considering the subsequence $n_k:= \lfloor k^{2/\varepsilon}\rfloor$, we see that $\sum_k P_\omega ( \tau_{n_k} \ge m(n_k) )< \infty$, $\P$-a.s. By the Borel--Cantelli lemma, for $\P$-almost all $\omega$ and $P_\omega$-almost all sufficiently large $k$, $\tau_{n_k} < m(n_k)$, which implies that for $n\in [n_{k-1}, n_k]$ and large $k$, we have $\tau_n < m(n_k) \le {n_k^{1+2\varepsilon} \over \E [\beta_{n_k} (e_1)]} \le {n^{1+3\varepsilon} \over \E [\beta_n(e_1)]}$ (the last inequality following from the estimate of $\E [\beta_n(e_1)]$ in Proposition \ref{p:beta-gamma}). In view of Proposition \ref{p:beta-gamma}, and since $\varepsilon$ can be as small as possible, this gives the lower bound in (\ref{nullrec}) of Theorem \ref{t:nullrec}.

To prove the upper bound, we note that $\alpha_{n,\lambda}(x) \le \beta_n(x)$ for any $\lambda\ge 0$ and any $0<|x|\le n$ (this is easily checked by induction on $|x|$). Thus, by (\ref{Laplace-tau}), for any $\lambda\ge 0$,
$$
E_\omega\left( \ee^{- \lambda \tau_n} \right) \le {\sum_{i=1}^\deg \omega (e, e_i) \beta_n (e_i) \over \sum_{i=1}^\deg \omega (e, e_i) \beta_{n,\lambda} (e_i)} \le \sum_{i=1}^\deg {\beta_n (e_i) \over \beta_{n,\lambda} (e_i)}.
$$

We now fix $r\in (1, \, {1\over \nu})$, where $\nu:= 1- {1\over \min\{ \kappa, \, 2\} }$ is defined in (\ref{theta}). It is possible to choose a small $\varepsilon>0$ such that
$$
{1\over \kappa -1} - {r\over \kappa}> 3\varepsilon \quad \hbox{if }\kappa \in (1, \, 2], \qquad 1 - {r\over 2}> 3\varepsilon \quad \hbox{if }\kappa \in (2, \, \infty].
$$

\noindent Let $\lambda = \lambda(n) := n^{-r}$. By (\ref{toto3}), we have $\beta_{n,n^{-r}} (e_i) \ge n^{-\varepsilon}\, \E [\beta_{n,n^{-r}} (e_i)]$ for $\P$-almost all $\omega$ and all sufficiently large $n$, which yields
$$
E_\omega\left( \ee^{- n^{-r} \tau_n} \right) \le n^\varepsilon \sum_{i=1}^\deg {\beta_n (e_i) \over \E [\beta_{n, n^{-r}} (e_i)]} .
$$

\noindent It is easy to bound $\beta_n (e_i)$. For any given $x\in \T \backslash \{ e\}$ with $|x|\le n$, $n\mapsto \beta_n (x)$ is non-increasing (this is easily checked by induction on $|x|$). Chebyshev's inequality, together with the Borel--Cantelli lemma (applied to a subsequence, as we did in the proof of the lower bound) and the monotonicity of $n\mapsto \beta_n(e_i)$, readily yields $\beta_n (e_i) \le n^\varepsilon \, \E [\beta_n (e_i)]$ for almost all $\omega$ and all sufficiently large $n$. As a consequence, for $\P$-almost all $\omega$ and all sufficiently large $n$,
$$
E_\omega\left( \ee^{- n^{-r} \tau_n} \right) \le n^{2\varepsilon} \sum_{i=1}^\deg {\E [\beta_n (e_i)] \over \E [\beta_{n, n^{-r}} (e_i)]} .
$$

\noindent By Proposition \ref{p:beta-gamma}, this yields $E_\omega ( \ee^{- n^{-r} \tau_n} ) \le n^{-\varepsilon}$ (for $\P$-almost all $\omega$ and all sufficiently large $n$; this is where we use ${1\over \kappa -1} - {r\over \kappa}> 3\varepsilon$ if $\kappa \in (1, \, 2]$, and $1 - {r\over 2}> 3\varepsilon$ if $\kappa \in (2, \, \infty]$). In particular, for $n_k:= \lfloor k^{2/\varepsilon} \rfloor$, we have $\P$-almost surely, $E_\omega ( \sum_k \ee^{- n_k^{-r} \tau_{n_k}} ) < \infty$, which implies that, $\p$-almost surely for all sufficiently large $k$, $\tau_{n_k} \ge n_k^r$. This implies that $\p$-almost surely for all sufficiently large $n$, $\tau_n \ge {1\over 2}\, n^r$. The upper bound in (\ref{nullrec}) of Theorem \ref{t:nullrec} follows.\hfill$\Box$

\bigskip

Proposition \ref{p:beta-gamma} is proved in Section \ref{s:beta-gamma}.

\section{Proof of Proposition \ref{p:beta-gamma}}
\label{s:beta-gamma}

Let $\theta \in [0,\, 1]$. Let $(Z_{n,\theta})$ be a sequence of random variables, such that $Z_{1,\theta} \; \buildrel law \over = \; \sum_{i=1}^\deg A_i$, where $(A_i, \, 1\le i\le \deg)$ is distributed as $(A(x_i), \, 1\le i\le \deg)$ (for any $x\in \T$), and that
\begin{equation}
    Z_{j+1,\theta} \; \buildrel law \over = \;
    \sum_{i=1}^\deg A_i {\theta +
    Z_{j,\theta}^{(i)}   \over
    1+ Z_{j,\theta}^{(i)} } , \qquad \forall\,
    j\ge 1,
    \label{ZW}
\end{equation}

\noindent where $Z_{j,\theta}^{(i)}$ (for $1\le i \le \deg$) are independent copies of $Z_{j,\theta}$, and are independent of the random vector $(A_i, \, 1\le i\le \deg)$.

Then, for any given $n\ge 1$ and $\lambda\ge 0$,
\begin{equation}
    Z_{n, 1-\ee^{-2\lambda}} \; \buildrel law \over
    = \;  \sum_{i=1}^\deg A_i\, \beta_{n,
    \lambda}(e_i) ,
    \label{Z=beta}
\end{equation}

\noindent provided $(A_i, \, 1\le i\le \deg)$ and $(\beta_{n, \lambda}(e_i), \, 1\le i\le \deg)$ are independent.

\medskip

\begin{proposition}
 \label{p:concentration}
 Assume $p={1\over \deg}$ and $\psi'(1)<0$. Let
 $\kappa$ be as in $(\ref{kappa})$. For all
 $a\in (1, \kappa) \cap (1, 2]$, we have
 $$
 \sup_{\theta \in [0,1]} \sup_{j\ge 1}
 {[\E (Z_{j,\theta} )^a ] \over (\E
 Z_{j,\theta})^a} < \infty.
 $$

\end{proposition}

\medskip

\noindent {\it Proof of Proposition \ref{p:concentration}.} Let $a\in (1,2]$. Conditioning on $A_1$, $\dots$, $A_\deg$, we can apply Lemma \ref{l:moment} to see that
\begin{eqnarray*}
 &&\E \left[ \left( \, \sum_{i=1}^\deg A_i
    {\theta+ Z_{j,\theta}^{(i)} \over 1+
    Z_{j,\theta}^{(i)} }
    \right)^a \Big| A_1, \dots, A_\deg \right]
    \\
 &\le& \sum_{i=1}^\deg A_i^a \, \E \left[ \left( {\theta+ Z_{j,\theta}
    \over 1+ Z_{j,\theta} }\right)^a \;
    \right] + (\deg-1) \left[ \sum_{i=1}^\deg A_i\,
    \E \left( {\theta+ Z_{j,\theta} \over 1+
    Z_{j,\theta} }
    \right) \right]^a
    \\
 &\le& \sum_{i=1}^\deg A_i^a \, \E \left[ \left( {\theta+ Z_{j,\theta}
    \over 1+ Z_{j,\theta} }\right)^a \;
    \right] + c_{32} \left[ \E \left( {\theta+ Z_{j,\theta} \over 1+
    Z_{j,\theta} } \right) \right]^a,
\end{eqnarray*}

\noindent where $c_{32}$ depends on $a$, $\deg$ and the bound of $A$ (recalling that $A$ is bounded away from 0 and infinity). Taking expectation on both sides, and in view of (\ref{ZW}), we obtain:
$$
\E[(Z_{j+1,\theta})^a] \le \deg \E(A^a) \E \left[ \left( {\theta+ Z_{j,\theta}\over
1+ Z_{j,\theta} }\right)^a \; \right] + c_{32} \left[ \E \left( {\theta+
Z_{j,\theta} \over 1+ Z_{j,\theta} } \right) \right]^a.
$$

\noindent We divide by $(\E Z_{j+1,\theta})^a = [ \E({\theta+Z_{j,\theta}\over 1+
Z_{j,\theta} })]^a$ on both sides, to see that
$$
{\E[(Z_{j+1,\theta})^a]\over (\E Z_{j+1,\theta})^a} \le \deg \E(A^a) {\E[ ({\theta+ Z_{j,\theta}
\over 1+ Z_{j,\theta} })^a] \over [\E ({\theta+ Z_{j,\theta} \over 1+ Z_{j,\theta} })]^a }
+ c_{32}.
$$

\noindent Put $\xi = \theta+ Z_{j,\theta}$. By (\ref{RSD}), we have
$$
{\E[ ({\theta+Z_{j,\theta} \over 1+Z_{j,\theta} })^a] \over [\E ({\theta+Z_{j,\theta} \over 1+ Z_{j,\theta} })]^a } = {\E[ ({\xi \over 1- \theta+ \xi })^a] \over [\E ({ \xi \over 1- \theta+ \xi })]^a } \le {\E[\xi^a] \over [\E \xi ]^a } .
$$

\noindent Applying Lemma \ref{l:moment} to $k=2$ yields that  $\E[\xi^a] =  \E[( \theta+ Z_{j,\theta} )^a] \le \theta^a + \E[( Z_{j,\theta} )^a] + (\theta + \E( Z_{j,\theta} ))^a $. It follows that ${\E[ \xi^a] \over [\E \xi ]^a }  \le  {\E[ (Z_{j,\theta})^a] \over [\E Z_{j,\theta}]^a } +2$, which implies that for $j\ge 1$,
$$
{\E[(Z_{j+1,\theta})^a]\over (\E Z_{j+1,\theta})^a} \le \deg \E(A^a) {\E[(Z_{j,\theta})^a]\over (\E Z_{j,\theta})^a} + (2 \deg \E(A^a)+ c_{32}).
$$

\noindent Thus, if $\deg \E(A^a)<1$ (which is the case if $1<a<\kappa$), then
$$
\sup_{j\ge 1} {\E[ (Z_{j,\theta})^a] \over (\E Z_{j,\theta})^a} < \infty,
$$

\noindent uniformly in $\theta \in [0, \, 1]$.\hfill$\Box$

\bigskip

We now turn to the proof of Proposition \ref{p:beta-gamma}. For the sake of clarity, the proofs of (\ref{comp-Laplace}), (\ref{E(beta):kappa>2}) and (\ref{E(beta):kappa<2}) are presented in three distinct parts.

\subsection{Proof of (\ref{comp-Laplace})}
\label{subs:beta}

By (\ref{exp}) and (\ref{ZW}), we have, for all $\theta\in [0, \, 1]$ and $j\ge 1$,
$$
\E \left\{ \exp\left( - t \, { Z_{j+1, \theta} \over \E (Z_{j+1, \theta})}\right) \right\} \le \E \left\{ \exp\left( - t  \sum_{i=1}^\deg A_i { Z^{(i)}_{j, \theta} \over \E (Z^{(i)}_{j, \theta}) }\right) \right\}, \qquad t\ge 0.
$$

\noindent Let $f_j(t) := \E \{ \exp ( - t { Z_{j, \theta} \over \E  Z_{j, \theta}} )\}$ and $g_j(t):= \E (\ee^{ -t\, M_j})$ (for $j\ge 1$). We have
$$
f_{j+1}(t) \le \E \left( \prod_{i=1}^\deg f_j(t A_i) \right), \quad j\ge 1.
$$

\noindent On the other hand, by (\ref{defMei1}),
$$
g_{j+1}(t) = \E \left\{ \exp\left( - t  \sum_{i=1}^\deg A(e_i) M^{(e_i)}_{j+1} \right) \right\} = \E \left( \prod_{i=1}^\deg g_j(t A_i) \right), \qquad j\ge 1.
$$

\noindent Since $f_1(\cdot)= g_1(\cdot)$, it follows by induction on $j$ that for all $j\ge 1$, $f_j(t) \le g_j(t)$; in particular, $f_n(t) \le g_n(t)$. We take $\theta = 1- \ee^{-2\lambda}$. In view of (\ref{Z=beta}), we have proved that
\begin{equation}
    \E \left\{ \exp\left( - t \sum_{i=1}^\deg A(e_i) 
    {\beta_{n, \lambda}(e_i) \over \E [\beta_{n, 
    \lambda}(e_i)] }\right) \right\} \le \E \left\{ 
    \ee^{- t \, M_n} \right\} ,
    \label{beta_n(e)}
\end{equation}

\noindent which yields (\ref{comp-Laplace}).\hfill$\Box$

\bigskip

\noindent {\bf Remark.} Let
$$
\beta_{n,\lambda}(e) := {(1-\ee^{-2\lambda})+ \sum_{i=1}^\deg A(e_i) \beta_{n,\lambda}(e_i) \over 1+ \sum_{i=1}^\deg A(e_i) \beta_{n,\lambda}(e_i)}.
$$

\noindent By (\ref{beta_n(e)}) and (\ref{exp}), if $\E(A)= {1\over \deg}$, then for $\lambda\ge 0$, $n\ge 1$ and $t\ge 0$, 
$$
\E \left\{ \exp\left( - t {\beta_{n, \lambda}(e) \over \E [\beta_{n, \lambda}(e)] }\right) \right\} \le \E \left\{ \ee^{- t \, M_n} \right\} .
$$

\subsection{Proof of (\ref{E(beta):kappa>2})}
\label{subs:kappa>2}

Assume $p={1\over \deg}$ and $\psi'(1)<0$. Since $Z_{j, \theta}$ is bounded uniformly in $j$, we have, by (\ref{ZW}), for $1\le j \le n-1$,
\begin{eqnarray}
    \E(Z_{j+1, \theta})
 &=& \E\left( {\theta+Z_{j, \theta} \over 1+Z_{j,
    \theta} } \right)
    \nonumber
    \\
 &\le& \E\left[(\theta+ Z_{j, \theta} )(1 - c_{33}\,
    Z_{j, \theta} )\right]
    \nonumber
    \\
 &\le & \theta + \E(Z_{j, \theta}) - c_{33}\,
    \E\left[(Z_{j, \theta})^2\right]
    \label{E(Z2)}
    \\
 &\le & \theta + \E(Z_{j, \theta}) - c_{33}\,
    \left[ \E Z_{j, \theta} \right]^2.
    \nonumber
\end{eqnarray}

\noindent By Lemma \ref{l:abc}, we have, for any $K>0$ and uniformly in $\theta\in [0, \,Ê{K\over n}]$,
\begin{equation}
    \label{53}
    \E (Z_{n, \theta}) \le c_{34} \left( \sqrt
    {\theta} + {1\over n} \right) \le {c_{35} \over
    \sqrt{n}}.
\end{equation}

\noindent We mention that this holds for all $\kappa \in (1, \, \infty]$. In view of (\ref{Z=beta}), this yields the upper bound in (\ref{E(beta):kappa>2}).

To prove the lower bound, we observe that
\begin{equation}
    \E(Z_{j+1, \theta}) \ge \E\left[(\theta+ Z_{j,
    \theta} )(1 - Z_{j, \theta} )\right] = \theta+
    (1-\theta) \E(Z_{j, \theta}) - \E\left[(Z_{j,
    \theta})^2\right] .
    \label{51}
\end{equation}

\noindent If furthermore $\kappa \in (2, \infty]$, then $\E [(Z_{j, \theta})^2 ] \le c_{36}\, (\E Z_{j, \theta})^2$ (see Proposition \ref{p:concentration}). Thus, for all $1\le j\le n-1$,
$$
\E(Z_{j+1, \theta}) \ge \theta+ (1-\theta) \E(Z_{j, \theta}) - c_{36}\, (\E Z_{j,\theta})^2 .
$$

\noindent By (\ref{53}), $\E (Z_{n, \theta}) \to 0$ uniformly in $\theta\in [0, \,Ê{K\over n}]$ (for any given $K>0$). An application of (\ref{Z=beta}) and Lemma \ref{l:abc} readily yields the lower bound in (\ref{E(beta):kappa>2}).\hfill$\Box$

\subsection{Proof of (\ref{E(beta):kappa<2})}
\label{subs:kappa<2}

We assume in this part $p={1\over \deg}$, $\psi'(1)<0$ and $1<\kappa \le 2$.

Let $\varepsilon>0$ be small. Since  $(Z_{j, \theta})$ is bounded, we have $\E[(Z_{j, \theta})^2] \le c_{37} \, \E [(Z_{j, \theta})^{\kappa-\varepsilon}]$, which, by Proposition \ref{p:concentration}, implies
\begin{equation}
    \E\left[ (Z_{j, \theta})^2 \right] \le c_{38} \,
    \left( \E Z_{j, \theta} \right)^{\kappa-
    \varepsilon} .
    \label{c38}
\end{equation}

\noindent Therefore, (\ref{51}) yields that
$$
\E(Z_{j+1, \theta}) \ge \theta+ (1-\theta) \E(Z_{j, \theta}) - c_{38} \, (\E Z_{j, \theta})^{\kappa-\varepsilon} .
$$

\noindent By (\ref{53}), $\E (Z_{n, \theta}) \to 0$ uniformly in $\theta\in [0, \,Ê{K\over n}]$ (for any given $K>0$). An application of Lemma \ref{l:abc} implies that for any $K>0$,
\begin{equation}
    \E (Z_{\ell, \theta}) \ge c_{14} \left(
    \theta^{1/(\kappa-\varepsilon)} + {1\over
    \ell^{1/(\kappa -1 - \varepsilon)}} \right),
    \qquad \forall \, \theta\in [0, \,Ê{K\over n}],
    \; \; \forall \, 1\le \ell \le n.
    \label{ell}
\end{equation}

\noindent The lower bound in (\ref{E(beta):kappa<2}) follows from (\ref{Z=beta}).

It remains to prove the upper bound. Define
$$
Y_{j, \theta} := {Z_{j, \theta} \over \E(Z_{j, \theta})} , \qquad 1\le j\le n.
$$

\noindent We take $Z_{j-1, \theta}^{(x)}$ (for $x\in \T_1$) to be independent copies of $Z_{j-1, \theta}$, and independent of $(A(x), \; x\in \T_1)$. By (\ref{ZW}), for $2\le j\le n$,
\begin{eqnarray*}
    Y_{j, \theta}
 &\; {\buildrel law \over =} \;& \sum_{x\in \T_1}
    A(x) {(\theta+ Z_{j-1, \theta}^{(x)} )/ (1+
    Z_{j-1, \theta}^{(x)}) \over \E
    [(\theta+ Z_{j-1, \theta}^{(x)} )/ (1+ Z_{j-1,
    \theta}^{(x)}) ]} \ge
    \sum_{x\in \T_1}
    A(x) {Z_{j-1, \theta}^{(x)} /
    (1+ Z_{j-1, \theta}^{(x)}) \over \theta+ \E
    [Z_{j-1, \theta}]}
    \\
 &=& { \E [Z_{j-1, \theta}]\over \theta+ \E
    [Z_{j-1, \theta}]} \sum_{x\in \T_1}
    A(x)
    Y_{j-1, \theta}^{(x)} - { \E [Z_{j-1,
    \theta}]\over \theta+ \E
    [Z_{j-1, \theta}]} \sum_{x\in \T_1}
    A(x) {(Z_{j-1, \theta}^{(x)})^2/\E(Z_{j-1,
    \theta}) \over 1+Z_{j-1, \theta}^{(x)}}
    \\
 &\ge& \sum_{x\in \T_1}
    A(x) Y_{j-1, \theta}^{(x)} -
    \Delta_{j-1, \theta} \; ,
\end{eqnarray*}

\noindent where
\begin{eqnarray*}
    Y_{j-1, \theta}^{(x)}
 &:=&{Z_{j-1, \theta}^{(x)} \over \E(Z_{j-1,
    \theta})} ,
    \\
    \Delta_{j-1, \theta}
 &:=&{\theta\over \theta+ \E [Z_{j-1, \theta}]}
    \sum_{x\in \T_1} A(x) Y_{j-1, \theta}^{(x)} +
    \sum_{x\in \T_1} A(x) {(Z_{j-1, \theta}^{(x)})^2
    \over \E(Z_{j-1, \theta})} .
\end{eqnarray*}

\noindent By (\ref{c38}), $\E[ {(Z_{j-1, \theta}^{(i)})^2 \over \E(Z_{j-1, \theta})}]\le c_{38}\, (\E Z_{j-1, \theta})^{\kappa-1-\varepsilon}$. On the other hand, by (\ref{ell}), $\E(Z_{j-1, \theta}) \ge c_{14}\, \theta^{1/(\kappa-\varepsilon)}$ for $2\le j \le n$, and thus ${\theta\over \theta+ \E [Z_{j-1, \theta}]} \le c_{39}\, (\E Z_{j-1, \theta})^{\kappa-1- \varepsilon}$. As a consequence, $\E( \Delta_{j-1, \theta} ) \le c_{40}\, (\E Z_{j-1, \theta})^{\kappa-1-\varepsilon}$.

If we write $\xi \; {\buildrel st. \over \ge} \; \eta$ to denote that $\xi$ is stochastically greater than or equal to $\eta$, then we have proved that $Y_{j, \theta} \; {\buildrel st. \over \ge} \; \sum_{x\in \T_1}^\deg A(x) Y_{j-1, \theta}^{(x)} - \Delta_{j-1, \theta}$. Applying the same argument to each of $(Y_{j-1, \theta}^{(x)}, \, x\in \T_1)$, we see that, for $3\le j\le n$,
$$
Y_{j, \theta} \; {\buildrel st. \over \ge} \; \sum_{u\in \T_1} A(u) \sum_{v\in \T_2: \; u={\buildrel \leftarrow \over v}} A(v) Y_{j-2, \theta}^{(v)} - \left( \Delta_{j-1, \theta}+ \sum_{u\in \T_1} A(u) \Delta_{j-2, \theta}^{(u)} \right) ,
$$

\noindent where $Y_{j-2, \theta}^{(v)}$ (for $v\in \T_2$) are independent copies of $Y_{j-2, \theta}$, and are independent of $(A(w), \, w\in \T_1 \cup \T_2)$, and $(\Delta_{j-2, \theta}^{(u)}, \, u\in \T_1)$ are independent of $(A(u), \, u\in \T_1)$ and are such that $\e[\Delta_{j-2, \theta}^{(u)}] \le c_{40}\, (\E Z_{j-2, \theta})^{\kappa-1-\varepsilon}$.

By induction, we arrive at: for $j>m \ge 1$,
\begin{equation}
    Y_{j, \theta} \; {\buildrel st. \over \ge}\;
    \sum_{x\in \T_m} \left( \prod_{y\in ]\! ] e, x
    ]\! ]} A(y) \right) Y_{j-m, \theta}^{(x)} -
    \Lambda_{j,m,\theta},
    \label{Yn>}
\end{equation}

\noindent where $Y_{j-m, \theta}^{(x)}$ (for $x\in \T_m$) are independent copies of $Y_{j-m, \theta}$, and are independent of the random vector $(A(w), \, 1\le |w| \le m)$, and $\E(\Lambda_{j,m,\theta}) \le c_{40}\, \sum_{\ell=1}^m (\E Z_{j-\ell, \theta})^{\kappa-1-\varepsilon} $.

Since $\E(Z_{i, \theta}) = \E({\theta+ Z_{i-1, \theta} \over 1+ Z_{i-1, \theta}}) \ge
\E(Z_{i-1, \theta}) - \E[(Z_{i-1, \theta})^2] \ge \E(Z_{i-1, \theta}) - c_{38}\, [\E Z_{i-1, \theta}
]^{\kappa-\varepsilon}$ (by (\ref{c38})), we have, for all $j\in (j_0, n]$ (with a large but fixed integer $j_0$) and $1\le \ell \le j-j_0$,
\begin{eqnarray*}
    \E(Z_{j, \theta})
 &\ge&\E(Z_{j-\ell, \theta}) \prod_{i=1}^\ell
    \left\{ 1- c_{38}\, [\E Z_{j-i, \theta}
    ]^{\kappa-1-\varepsilon}\right\}
    \\
 &\ge&\E(Z_{j-\ell, \theta}) \prod_{i=1}^\ell
    \left\{ 1- c_{41}\, (j-i)^{-(\kappa-1-
    \varepsilon)/2}\right\} ,
\end{eqnarray*}

\noindent the last inequality being a consequence of (\ref{53}). Thus, for $j\in (j_0, n]$ and $1\le \ell \le j^{(\kappa-1-\varepsilon)/2}$, $\E(Z_{j, \theta}) \ge c_{42}\, \E(Z_{j-\ell, \theta})$, which implies that for all $m\le j^{(\kappa-1-\varepsilon)/2}$, $\E(\Lambda_{j,m, \theta}) \le c_{43} \, m (\E Z_{j, \theta})^{\kappa-1-\varepsilon}$. By Chebyshev's inequality, for $j\in (j_0, n]$, $m\le j^{(\kappa-1-\varepsilon)/2}$ and  $r>0$,
\begin{equation}
    \P\left\{ \Lambda_{j,m, \theta} > \varepsilon
    r\right\} \le {c_{43} \, m (\E Z_{j,
    \theta})^{\kappa -1-\varepsilon} \over
    \varepsilon r}.
    \label{toto4}
\end{equation}

Let us go back to (\ref{Yn>}), and study the behaviour of $\sum_{x\in \T_m} ( \prod_{y\in ]\! ] e, x ]\! ]} A(y) ) Y_{j-m, \theta}^{(x)}$. Let $M^{(x)}$ (for $x\in \T_m$) be independent copies of $M_\infty$ and independent of all other random variables. Since $\E(Y_{j-m, \theta}^{(x)})= \E(M^{(x)})=1$, we have, by Fact \ref{f:petrov}, for any $a\in (1, \, \kappa)$,
\begin{eqnarray*}
 &&\E \left\{ \left| \sum_{x\in \T_m} \left(
    \prod_{y\in ]\! ] e, x ]\! ]} A(y) \right)
    (Y_{j- m, \theta}^{(x)} - M^{(x)}) \right|^a
    \right\}
    \\
 &\le&2 \E \left\{ \sum_{x\in \T_m} \left(
    \prod_{y\in ]\! ] e, x ]\! ]} A(y)^a \right)
    \, \E\left( | Y_{j-m, \theta}^{(x)} -
    M^{(x)}|^a \right) \right\}.
\end{eqnarray*}

\noindent By Proposition \ref{p:concentration} and the fact that $(M_n)$ is a martingale bounded in $L^a$, we have $\E ( | Y_{j-m, \theta}^{(x)} - M^{(x)}|^a ) \le c_{44}$. Thus,
\begin{eqnarray*}
    \E \left\{ \left| \sum_{x\in \T_m} \left(
    \prod_{y\in ]\! ] e, x ]\! ]} A(y) \right)
    (Y_{j- m, \theta}^{(x)} - M^{(x)}) \right|^a
    \right\}
 &\le& 2c_{44} \E \left\{ \sum_{x\in \T_m}
    \prod_{y\in ]\! ] e, x ]\! ]} A(y)^a \right\}
    \\
 &=& 2c_{44} \, \deg^m \, [\E(A^a)]^m.
\end{eqnarray*}

\noindent By Chebyshev's inequality,
\begin{equation}
    \P \left\{ \left| \sum_{x\in \T_m} \left(
    \prod_{y\in ]\! ] e, x ]\! ]} A(y) \right)
    (Y_{j- m, \theta}^{(x)} - M^{(x)}) \right| >
    \varepsilon r\right\} \le {2c_{44} \, \deg^m
    [\E(A^a)]^m \over \varepsilon^a r^a}.
    \label{toto6}
\end{equation}

\noindent Clearly, $\sum_{x\in \T_m} (\prod_{y\in ]\! ] e, x ]\! ]} A(y) ) M^{(x)}$ is distributed as $M_\infty$. We can thus plug (\ref{toto6}) and (\ref{toto4}) into (\ref{Yn>}), to see that for $j\in [j_0, n]$, $m\le j^{(\kappa-1-\varepsilon)/2}$ and $r>0$,
\begin{equation}
    \P \left\{ Y_{j, \theta} > (1-2\varepsilon)
    r\right\} \ge \P \left\{ M_\infty > r\right\} -
    {c_{43}\, m (\E Z_{j, \theta})^{\kappa-1-
    \varepsilon} \over \varepsilon r} - {2c_{44} \,
    \deg^m [\E(A^a)]^m \over
    \varepsilon^a r^a} .
    \label{Yn-lb}
\end{equation}

\noindent We choose $m:= \lfloor j^\varepsilon \rfloor$. Since $a\in (1, \, \kappa)$, we have $\deg \E(A^a) <1$, so that $\deg^m [\E(A^a)]^m \le \exp( - j^{\varepsilon/2})$ for all large $j$. We choose $r= {1\over (\E Z_{j, \theta})^{1- \delta}}$, with $\delta := {4\kappa \varepsilon \over \kappa -1}$. In view of (\ref{M-tail}), we obtain: for $j\in [j_0, n]$,
$$
\P \left\{ Y_{j, \theta} > {1-2\varepsilon\over (\E Z_{j, \theta})^{1- \delta}} \right\} \ge c_{23} \, (\E Z_{j, \theta})^{(1- \delta) \kappa} - {c_{43}\over \varepsilon} \, j^\varepsilon\, (\E Z_{j, \theta})^{\kappa-\varepsilon-\delta} - {2c_{44} \, (\E Z_{j, \theta})^{(1- \delta)a} \over \varepsilon^a \exp(j^{\varepsilon/2})} .
$$

\noindent Since $c_{14}/j^{1/(\kappa-1- \varepsilon)} \le \E(Z_{j, \theta}) \le c_{35}/j^{1/2}$ (see (\ref{ell}) and (\ref{53}), respectively), we can pick up sufficiently small $\varepsilon$, so that for $j\in [j_0, n]$,
$$
\P \left\{ Y_{j, \theta} > {1-2\varepsilon\over (\E Z_{j, \theta})^{1- \delta}} \right\} \ge {c_{23} \over 2} \, (\E Z_{j, \theta})^{(1-\delta) \kappa}.
$$

\noindent Recall that by definition, $Y_{j, \theta} = {Z_{j, \theta} \over \E(Z_{j, \theta})}$. Therefore, for $j\in [j_0, n]$,
$$
\E[(Z_{j, \theta})^2] \ge [\E Z_{j, \theta}]^2 \, {(1-2\varepsilon)^2\over  (\E Z_{j, \theta})^{2(1- \delta)}} \P \left\{ Y_{j, \theta} > {1-2\varepsilon \over (\E Z_{j, \theta})^{1- \delta}} \right\} \ge c_{45} \, (\E Z_{j, \theta})^{\kappa+ (2- \kappa)\delta}.
$$

\noindent Of course, the inequality holds trivially for $0\le j < j_0$ (with possibly a different value of the constant $c_{45}$). Plugging this into (\ref{E(Z2)}), we see that for $1\le j\le n-1$,
$$
\E(Z_{j+1, \theta}) \le \theta + \E(Z_{j, \theta}) - c_{46}\, (\E Z_{j, \theta})^{\kappa+ (2- \kappa)\delta} .
$$

\noindent By Lemma \ref{l:abc}, this yields $\E(Z_{n, \theta}) \le c_{47} \, \{ \theta^{1/[\kappa+ (2- \kappa)\delta]} + n^{- 1/ [\kappa -1 + (2- \kappa)\delta]}\}$. An application of (\ref{Z=beta}) implies the desired upper bound in (\ref{E(beta):kappa<2}).\hfill$\Box$

\bigskip

\noindent {\bf Remark.} A close inspection on our argument shows that under the assumptions $p= {1\over \deg}$ and $\psi'(1)<0$, we have, for any $1\le i \le \deg$ and uniformly in $\lambda \in [0, \, {1\over n}]$,
$$
\left( {\alpha_{n, \lambda}(e_i) \over \E[\alpha_{n, \lambda}(e_i)]} ,\; {\beta_{n, \lambda}(e_i) \over \E[\beta_{n, \lambda}(e_i)]} , \; {\gamma_n(e_i) \over \E[\gamma_n (e_i)]} \right) \; {\buildrel law \over \longrightarrow} \; (M_\infty, \, M_\infty, \, M_\infty),
$$

\noindent where ``${\buildrel law \over \longrightarrow}$" stands for convergence in distribution, and $M_\infty$ is the random variable defined in $(\ref{cvg-M})$.\hfill$\Box$

\bigskip
\bigskip

\noindent {\Large\bf Acknowledgements}

\bigskip

We are grateful to Philippe Carmona and Marc Yor for helpful discussions.


\bigskip
\bigskip


{\footnotesize

\baselineskip=12pt

\noindent
\begin{tabular}{lll}
&\hskip20pt Yueyun Hu
    & \hskip45pt Zhan Shi\\
&\hskip20pt D\'epartement de Math\'ematiques
    & \hskip45pt Laboratoire de Probabilit\'es et
       Mod\`eles Al\'eatoires \\
&\hskip20pt Universit\'e Paris XIII
    & \hskip45pt Universit\'e Paris VI\\
&\hskip20pt 99 avenue J-B Cl\'ement
    & \hskip45pt 4 place Jussieu\\
&\hskip20pt F-93430 Villetaneuse
    & \hskip45pt F-75252 Paris Cedex 05\\
&\hskip20pt France
    & \hskip45pt France \\
&\hskip20pt {\tt yueyun@math.univ-paris13.fr}
    & \hskip45pt {\tt zhan@proba.jussieu.fr}
\end{tabular}

}

\end{document}